\documentclass[final,3p]{elsarticle}
 \usepackage{graphics}
 \usepackage{graphicx}
 \usepackage{epsfig}
\usepackage{amssymb}
 \usepackage{amsthm}
 \usepackage{lineno}
 \usepackage{amsmath}
   \numberwithin{equation}{section}
\usepackage{mathrsfs}

\journal{Elsevier} 

\newtheorem{thm}{Theorem}[section]

\newtheorem{lem}[thm]{Lemma}
\newtheorem{prop}[thm]{Proposition}

\biboptions{sort&compress,square}
\allowdisplaybreaks
\begin{document}
\begin{frontmatter}
\author{Tong Wu$^{a}$}
\ead{wut977@nenu.edu.cn}
\author{Yong Wang$^{b,*}$}
\ead{wangy581@nenu.edu.cn}
\cortext[cor]{Corresponding author.}
\address{$^a$Department of Mathematics, Northeastern University, Shenyang, 110004, China}
\address{$^b$School of Mathematics and Statistics, Northeast Normal University,
Changchun, 130024, China}

\title{The Hodge-Dirac operator and Dabrowski-Sitarz-Zalecki type theorems for manifolds with boundary}
\begin{abstract}
In \cite{DL1}, Dabrowski etc. gave spectral Einstein bilinear functionals of differential forms for the Hodge-Dirac operator $d+\delta$ on an oriented even-dimensional Riemannian manifold. In this paper, we generalize the results of
Dabrowski etc. to the cases of $4$ dimensional oriented Riemannian manifolds with boundary. Furthermore, we give the proof of Dabrowski-Sitarz-Zalecki type theorems associated with the Hodge-Dirac operator for manifolds with boundary.
\end{abstract}
\end{frontmatter}
\section{Introduction}
The theory of noncommutative residue for one-dimensional manifolds was discovered by Manin \cite{Y} and
Adler \cite{A1} in connection with geometric aspects of nonlinear partial differential equations. For arbitrary
closed compact n-dimensional manifolds, the noncommutative reside was introduced by Wodzicki in \cite{Wo,Wo1},
using the theory of zeta functions of elliptic pseudodifferential operators. Let E be a finite-dimensional
complex vector bundle over a closed compact manifold M of dimension n, the noncommutative residue of a
pseudo-differential operator P$\in \Psi{DO}(E)$ can be defined by $${\rm Wres(P)}:=(2\pi)^{-n}\int_{S^*M}{\rm Tr}[\sigma^P_{-n}(x,\xi)]dxd\xi$$
where $S^*M \subset T^*M$ denotes the co-sphere bundle on $M$ and $\sigma^P_{-n}$ is the component of order $-n$ of the
complete symbol
$$\sigma^P:=\sum_i\sigma^P_i.$$
of P, and the linear functional ${\rm Wres} : \Psi{DO}(E) \rightarrow C$ is in fact the unique trace (up to multiplication by
constants) on the algebra of pseudo-differential operators $\Psi{DO}(E)$.
In \cite{Co1}, Connes computed a conformal four-dimensional Polyakov action analogy using the noncommutative residue. Connes proved that the noncommutative residue on a compact manifold $M$ coincided with
Dixmier's trace on pseudodifferential operators of order-dim$M$ in \cite{Co2,Co3}. The theory has very rich structures
both in physics and mathematics. More precisely, Connes made a challenging observation that the Wodzicki
residue of the inverse square of the Dirac operator yields the Einstein-Hilbert action of general relativity.
Kastler\cite{Ka} gave a brute-force proof of this theorem, and Kalau and Walze\cite{KW} proved this theorem in the
normal coordinates system simultaneously, which is called the Kastler-Kalau-Walze theorem now. Let $s$ be
the scalar curvature and ${\rm Wres}$ denotes the noncommutative residue, then the Kastler-Kalau-Walze theorem
gives an operator-theoretic explanation of the gravitational action and says that for a 4-dimensional closed
spin manifold and Dirac operator D, there exists a constant $c_0$, such that $${\rm Wres(D^{-2})}=c_0\int_Ms{\rm dVol_M}. $$

On the other hand, Fedosov etc. defined a noncommutative residue on Boutet de Monvel's algebra and
proved that it was a unique continuous trace in \cite{FGLS}, and generalized the definition of noncommutative
residue to manifolds with boundary. In \cite{ES}, Schrohe gave the relation between the Dixmier trace and the
noncommutative residue for manifolds with boundary. For elliptic pseudodifferential operators, Wang proved
the Kastler-Kalau-Walze type theorem and gave the operator-theoretic explanation of the gravitational
action for lower dimensional manifolds with boundary\cite{Wa1,Wa3,Wa4}.

In the noncommutative realm the spectral-theoretic approach to scalar curvature has been extended
also to quantum tori in the seminal work of Connes and Moscovici\cite{Co3}. Furthermore, the pseudodifferential
operators and symbol calculus introduced in \cite{Co4} and extended to crossed product algebras in \cite{B1,B2}, have
been employed for computations of certain values and residues of zeta functions of suitable Laplace type
operators. Recently, in order to recover other important tensors in both the classical setup as well as for the
generalised or quantum geometries, for the metric tensor g, Ricci curvature Ric and the scalar curvature s,
Dabrowski etc. \cite{DL} defined bilinear functionals $$G:=Ric-\frac{1}{2}s(g)g,$$
and they demonstrated that the noncommutative residue density recovered the tensors g and G as certain
bilinear functionals of vector fields on a manifold M, while their dual tensores are recovered as a density
of bilinear functionals of differential one-forms on M. Motivated by spectral Einstein bilinear functionals of differential forms for the Hodge-Dirac operator $d+\delta$ on an oriented even-dimensional Riemannian manifold in Proposition 3.3. \cite{DL1} and the Kastler-Kalau-Walze type theorem\cite{Ka,KW}, we give some
new spectral functionals which are the extension of spectral functionals for the Hodge-Dirac operator with Clifford
multiplication by the local coframe basis, and we relate them to the noncommutative residue for manifolds
with boundary. For lower dimensional compact Riemannian manifolds with boundary, we compute the
residue $\widetilde{{\rm Wres}}[\pi^+(c(w)(\widetilde{D}c(v)+c(v)\widetilde{D})\widetilde{D}^{-1})\circ\pi^+(\widetilde{D}^{-2})]$ and $\widetilde{{\rm Wres}}[\pi^+(c(w)(\widetilde{D}c(v)+c(v)\widetilde{D})\widetilde{D}^{-2})\circ\pi^+(\widetilde{D}^{-1})]$, which we call type-I operator and type-II operator  
and obtain the Dabrowski-Sitarz-Zalecki type
theorems for four dimensional oriented Riemannian manifolds with boundary. Our main theorems are as follows.
\begin{thm}\label{uo0o0o}
Let M be a $4$-dimensional compact oriented Riemannian manifold with boundary $\partial M$ and the metric
$g^M$ be defined in Section \ref{section:2}, then we get the following equality:
\begin{align}
\label{bo0o00}
&\widetilde{{\rm Wres}}[\pi^+(c(w)(\widetilde{D}c(v)+c(v)\widetilde{D})\widetilde{D}^{-1})\circ\pi^+(\widetilde{D}^{-2})]\nonumber\\
&=\frac{64\pi^2}{3}\int_M[Ric(v,w)-\frac{1}{2}s(g)g(v,w)]{Vol_g}+\int_{\partial M}\bigg\{-8\bigg(\sum_{j=1}^ng(e_j,\nabla^L_{e_j}v)g(w,\frac{\partial}{\partial{x_n}})-g(w,\nabla^L_{\frac{\partial}{\partial{x_n}}}v)\nonumber\\
&+g(\nabla^L_wv,\frac{\partial}{\partial{x_n}})\bigg) -\frac{8}{3}\partial_{x_n}g(v^T,w^T)+8\partial_{x_n}(v_nw_n)+\frac{88}{9}Kv_nw_n-\frac{36}{9}Kg(v^T,w^T)\nonumber\\
&+8<\nabla^L_v\frac{\partial}{\partial{x_n}},w^T>\bigg\}\pi^2dx',
\end{align}
where $c(v)$, $c(w)$, $v^T$ and $w^T$ are defined in Section \ref{section:3}.
\end{thm}
\begin{thm}\label{momo}
Let M be a $4$-dimensional compact oriented Riemannian manifold with boundary $\partial M$ and the metric
$g^M$ be defined in Section \ref{section:2}, then we get the following equality:
\begin{align}
\label{bmomo}
&\widetilde{{\rm Wres}}[\pi^+(c(w)(\widetilde{D}c(v)+c(v)\widetilde{D})\widetilde{D}^{-2})\circ\pi^+(\widetilde{D}^{-1})]\nonumber\\
&=\frac{64\pi^2}{3}\int_M[Ric(v,w)-\frac{1}{2}s(g)g(v,w)]{Vol_g}+\int_{\partial M}\bigg\{8\bigg(\sum_{j=1}^ng(e_j,\nabla^L_{e_j}v)g(w,\frac{\partial}{\partial{x_n}})-g(w,\nabla^L_{\frac{\partial}{\partial{x_n}}}v)\nonumber\\
&+g(\nabla^L_wv,\frac{\partial}{\partial{x_n}})\bigg) -\frac{8}{3}\partial_{x_n}g(v^T,w^T)-8\partial_{x_n}(v_nw_n)+\frac{32}{9}Kv_nw_n-\frac{16}{9}Kg(v^T,w^T)\nonumber\\
&-8<\nabla^L_v\frac{\partial}{\partial{x_n}},w^T>\bigg\}\pi^2dx',
\end{align}
where $c(v)$, $c(w)$, $v^T$ and $w^T$ are defined in Section \ref{section:3}.
\end{thm}
We note that our theorems may be generalied to general even dimensional manifolds, we plan to generalize our theorems to general even dimensional manifolds in the future.
The paper is organized in the following way. In Section \ref{section:2},  we recall some basic facts and formulas about Boutet de
Monvel's calculus and the definition of the noncommutative residue for manifolds with boundary. In Section \ref{section:3}, we recall the
 spectral Einstein bilinear functionals of differential forms for the Hodge-Dirac operator $d+\delta$ \cite{DL1}. In Section \ref{section:4}, we prove the Dabrowski-Sitarz-Zalecki type theorem associated with the residue $\widetilde{{\rm Wres}}[\pi^+(c(w)(\widetilde{D}c(v)+c(v)\widetilde{D})\widetilde{D}^{-1})\circ\pi^+(\widetilde{D}^{-2})]$ on manifolds with boundary. In Section \ref{section:5},  we prove the Dabrowski-Sitarz-Zalecki type theorem associated with the residue $\widetilde{{\rm Wres}}[\pi^+(c(w)(\widetilde{D}c(v)+c(v)\widetilde{D})\widetilde{D}^{-2})\circ\pi^+(\widetilde{D}^{-1})]$ on manifolds with boundary.
\section{Boutet de Monvel's calculus}
\label{section:2}
 In this section, we recall some basic facts and formulas about Boutet de
Monvel's calculus and the definition of the noncommutative residue for manifolds with boundary which will be used in the following. For more details, see Section 2 in \cite{Wa3}.\\
 \indent Let $M$ be a 4-dimensional compact oriented manifold with boundary $\partial M$.
We assume that the metric $g^{TM}$ on $M$ has the following form near the boundary,
\begin{equation}
\label{b1}
g^{M}=\frac{1}{h(x_{n})}g^{\partial M}+dx _{n}^{2},
\end{equation}
where $g^{\partial M}$ is the metric on $\partial M$ and $h(x_n)\in C^{\infty}([0, 1)):=\{\widehat{h}|_{[0,1)}|\widehat{h}\in C^{\infty}((-\varepsilon,1))\}$ for
some $\varepsilon>0$ and $h(x_n)$ satisfies $h(x_n)>0$, $h(0)=1$ where $x_n$ denotes the normal directional coordinate. Let $U\subset M$ be a collar neighborhood of $\partial M$ which is diffeomorphic with $\partial M\times [0,1)$. By the definition of $h(x_n)\in C^{\infty}([0,1))$
and $h(x_n)>0$, there exists $\widehat{h}\in C^{\infty}((-\varepsilon,1))$ such that $\widehat{h}|_{[0,1)}=h$ and $\widehat{h}>0$ for some
sufficiently small $\varepsilon>0$. Then there exists a metric $g'$ on $\widetilde{M}=M\bigcup_{\partial M}\partial M\times
(-\varepsilon,0]$ which has the form on $U\bigcup_{\partial M}\partial M\times (-\varepsilon,0 ]$
\begin{equation}
\label{b2}
g'=\frac{1}{\widehat{h}(x_{n})}g^{\partial M}+dx _{n}^{2} ,
\end{equation}
such that $g'|_{M}=g$. We fix a metric $g'$ on the $\widetilde{M}$ such that $g'|_{M}=g$.

Let the Fourier transformation $F'$  be
\begin{equation*}
F':L^2({\bf R}_t)\rightarrow L^2({\bf R}_v);~F'(u)(v)=\int_\mathbb{R} e^{-ivt}u(t)dt
\end{equation*}
and let
\begin{equation*}
r^{+}:C^\infty ({\bf R})\rightarrow C^\infty (\widetilde{{\bf R}^+});~ f\rightarrow f|\widetilde{{\bf R}^+};~
\widetilde{{\bf R}^+}=\{x\geq0;x\in {\bf R}\}.
\end{equation*}
\indent We define $H^+=F'(\Phi(\widetilde{{\bf R}^+}));~ H^-_0=F'(\Phi(\widetilde{{\bf R}^-}))$ which satisfies
$H^+\bot H^-_0$, where $\Phi(\widetilde{{\bf R}^+}) =r^+\Phi({\bf R})$, $\Phi(\widetilde{{\bf R}^-}) =r^-\Phi({\bf R})$ and $\Phi({\bf R})$
denotes the Schwartz space. We have the following
 property: $h\in H^+~$ (resp. $H^-_0$) if and only if $h\in C^\infty({\bf R})$ which has an analytic extension to the lower (resp. upper) complex
half-plane $\{{\rm Im}\xi<0\}$ (resp. $\{{\rm Im}\xi>0\})$ such that for all nonnegative integer $l$,
 \begin{equation*}
\frac{d^{l}h}{d\xi^l}(\xi)\sim\sum^{\infty}_{k=1}\frac{d^l}{d\xi^l}(\frac{c_k}{\xi^k}),
\end{equation*}
as $|\xi|\rightarrow +\infty,{\rm Im}\xi\leq0$ (resp. ${\rm Im}\xi\geq0)$ and where $c_k\in\mathbb{C}$ are some constants.\\
 \indent Let $H'$ be the space of all polynomials and $H^-=H^-_0\bigoplus H';~H=H^+\bigoplus H^-.$ Denote by $\pi^+$ (resp. $\pi^-$) the
 projection on $H^+$ (resp. $H^-$). Let $\widetilde H=\{$rational functions having no poles on the real axis$\}$. Then on $\tilde{H}$,
 \begin{equation}
 \label{b3}
\pi^+h(\xi_0)=\frac{1}{2\pi i}\lim_{u\rightarrow 0^{-}}\int_{\Gamma^+}\frac{h(\xi)}{\xi_0+iu-\xi}d\xi,
\end{equation}
where $\Gamma^+$ is a Jordan closed curve
included ${\rm Im}(\xi)>0$ surrounding all the singularities of $h$ in the upper half-plane and
$\xi_0\in {\bf R}$. In our computations, we only compute $\pi^+h$ for $h$ in $\widetilde{H}$. Similarly, define $\pi'$ on $\tilde{H}$,
\begin{equation}
\label{b4}
\pi'h=\frac{1}{2\pi}\int_{\Gamma^+}h(\xi)d\xi.
\end{equation}
So $\pi'(H^-)=0$. For $h\in H\bigcap L^1({\bf R})$, $\pi'h=\frac{1}{2\pi}\int_{{\bf R}}h(v)dv$ and for $h\in H^+\bigcap L^1({\bf R})$, $\pi'h=0$.\\
\indent An operator of order $m\in {\bf Z}$ and type $d$ is a matrix\\
$$\widetilde{A}=\left(\begin{array}{lcr}
  \pi^+P+G  & K  \\
   T  &  \widetilde{S}
\end{array}\right):
\begin{array}{cc}
\   C^{\infty}(M,E_1)\\
 \   \bigoplus\\
 \   C^{\infty}(\partial{M},F_1)
\end{array}
\longrightarrow
\begin{array}{cc}
\   C^{\infty}(M,E_2)\\
\   \bigoplus\\
 \   C^{\infty}(\partial{M},F_2)
\end{array},
$$
where $M$ is a manifold with boundary $\partial M$ and
$E_1,E_2$~ (resp. $F_1,F_2$) are vector bundles over $M~$ (resp. $\partial M
$).~Here,~$P:C^{\infty}_0(\Omega,\overline {E_1})\rightarrow
C^{\infty}(\Omega,\overline {E_2})$ is a classical
pseudodifferential operator of order $m$ on $\Omega$, where
$\Omega$ is a collar neighborhood of $M$ and
$\overline{E_i}|M=E_i~(i=1,2)$. $P$ has an extension:
$~{\cal{E'}}(\Omega,\overline {E_1})\rightarrow
{\cal{D'}}(\Omega,\overline {E_2})$, where
${\cal{E'}}(\Omega,\overline {E_1})~({\cal{D'}}(\Omega,\overline
{E_2}))$ is the dual space of $C^{\infty}(\Omega,\overline
{E_1})~(C^{\infty}_0(\Omega,\overline {E_2}))$. Let
$e^+:C^{\infty}(M,{E_1})\rightarrow{\cal{E'}}(\Omega,\overline
{E_1})$ denote extension by zero from $M$ to $\Omega$ and
$r^+:{\cal{D'}}(\Omega,\overline{E_2})\rightarrow
{\cal{D'}}(\Omega, {E_2})$ denote the restriction from $\Omega$ to
$X$, then define
$$\pi^+P=r^+Pe^+:C^{\infty}(M,{E_1})\rightarrow {\cal{D'}}(\Omega,
{E_2}).$$ In addition, $P$ is supposed to have the
transmission property; this means that, for all $j,k,\alpha$, the
homogeneous component $p_j$ of order $j$ in the asymptotic
expansion of the
symbol $p$ of $P$ in local coordinates near the boundary satisfies:\\
$$\partial^k_{x_n}\partial^\alpha_{\xi'}p_j(x',0,0,+1)=
(-1)^{j-|\alpha|}\partial^k_{x_n}\partial^\alpha_{\xi'}p_j(x',0,0,-1),$$
then $\pi^+P:C^{\infty}(M,{E_1})\rightarrow C^{\infty}(M,{E_2})$. Let $G$,$T$ be respectively the singular Green operator
and the trace operator of order $m$ and type $d$. Let $K$ be a
potential operator and $S$ be a classical pseudodifferential
operator of order $m$ along the boundary. Denote by $B^{m,d}$ the collection of all operators of
order $m$
and type $d$,  and $\mathcal{B}$ is the union over all $m$ and $d$.\\
\indent Recall that $B^{m,d}$ is a Fr\'{e}chet space. The composition
of the above operator matrices yields a continuous map:
$B^{m,d}\times B^{m',d'}\rightarrow B^{m+m',{\rm max}\{
m'+d,d'\}}.$ Write $$\widetilde{A}=\left(\begin{array}{lcr}
 \pi^+P+G  & K \\
 T  &  \widetilde{S}
\end{array}\right)
\in B^{m,d},
 \widetilde{A}'=\left(\begin{array}{lcr}
\pi^+P'+G'  & K'  \\
 T'  &  \widetilde{S}'
\end{array} \right)
\in B^{m',d'}.$$\\
 The composition $\widetilde{A}\widetilde{A}'$ is obtained by
multiplication of the matrices (For more details see \cite{ES}). For
example $\pi^+P\circ G'$ and $G\circ G'$ are singular Green
operators of type $d'$ and
$$\pi^+P\circ\pi^+P'=\pi^+(PP')+L(P,P').$$
Here $PP'$ is the usual
composition of pseudodifferential operators and $L(P,P')$ called
leftover term is a singular Green operator of type $m'+d$. For our case, $P,P'$ are classical pseudo differential operators, in other words $\pi^+P\in \mathcal{B}^{\infty}$ and $\pi^+P'\in \mathcal{B}^{\infty}$ .\\
\indent Let $M$ be a $n$-dimensional compact oriented manifold with boundary $\partial M$.
Denote by $\mathcal{B}$ the Boutet de Monvel's algebra. We recall that the main theorem in \cite{FGLS,Wa3}.
\begin{thm}\label{th:32}{\rm\cite{FGLS}}{\bf(Fedosov-Golse-Leichtnam-Schrohe)}
 Let $M$ and $\partial M$ be connected, ${\rm dim}M=n\geq3$, and let $\widetilde{S}$ (resp. $\widetilde{S}'$) be the unit sphere about $\xi$ (resp. $\xi'$) and $\sigma(\xi)$ (resp. $\sigma(\xi')$) be the corresponding canonical
$n-1$ (resp. $(n-2)$) volume form.
 Set $\widetilde{A}=\left(\begin{array}{lcr}\pi^+P+G &   K \\
T &  \widetilde{S}    \end{array}\right)$ $\in \mathcal{B}$ , and denote by $p$, $b$ and $s$ the local symbols of $P,G$ and $\widetilde{S}$ respectively.
 Define:
 \begin{align}
{\rm{\widetilde{Wres}}}(\widetilde{A})&=\int_X\int_{\bf \widetilde{ S}}{\rm{tr}}_E\left[p_{-n}(x,\xi)\right]\sigma(\xi)dx \nonumber\\
&+2\pi\int_ {\partial X}\int_{\bf \widetilde{S}'}\left\{{\rm tr}_E\left[({\rm{tr}}b_{-n})(x',\xi')\right]+{\rm{tr}}
_F\left[s_{1-n}(x',\xi')\right]\right\}\sigma(\xi')dx',
\end{align}
where ${\rm{\widetilde{Wres}}}$ denotes the noncommutative residue of an operator in the Boutet de Monvel's algebra.\\
Then~~ a) ${\rm \widetilde{Wres}}([\widetilde{A},B])=0 $, for any
$\widetilde{A},B\in\mathcal{B}$;~~ b) It is the unique continuous trace on
$\mathcal{B}/\mathcal{B}^{-\infty}$.
\end{thm}
\section{The spectral Einstein functional associated with the Hodge-Dirac operator} 
\label{section:3}
Firstly we recall the definition of the Hodge-Dirac operator. Let $M$ be an $n$-dimensional oriented compact Riemannian manifold with a Riemannian metric $g^{M}$ and let $\nabla^L$ be the Levi-Civita connection about $g^{M}$. In the fixed orthonormal frame $\{\widetilde{e}_1,\cdots,\widetilde{e}_n\}$, the connection matrix $(\omega_{s,t})$ is defined by
\begin{equation}
\label{a2}
\nabla^L(\widetilde{e}_1,\cdots,\widetilde{e}_n)= (\widetilde{e}_1,\cdots,\widetilde{e}_n)(\omega_{s,t}).
\end{equation}
\indent Let $\epsilon (\widetilde{e}_j^*)$,~$\iota (\widetilde{e}_j^*)$ be the exterior and interior multiplications respectively, where $\widetilde{e}_j^*=g^{TM}(\widetilde{e}_j,\cdot)$.
Write
\begin{equation}
\label{a3}
\widehat{c}(\widetilde{e}_j)=\epsilon (\widetilde{e}_j^* )+\iota
(e_j^*);~~
c(\widetilde{e}_j)=\epsilon (\widetilde{e}_j^* )-\iota (\widetilde{e}_j^* ),
\end{equation}
which satisfies
\begin{align}
\label{a4}
&\widehat{c}(\widetilde{e}_i)\widehat{c}(\widetilde{e}_j)+\widehat{c}(\widetilde{e}_j)\widehat{c}(\widetilde{e}_i)=2g^{M}(\widetilde{e}_i,\widetilde{e}_j);~~\nonumber\\
&c(\widetilde{e}_i)c(\widetilde{e}_j)+c(\widetilde{e}_j)c(\widetilde{e}_i)=-2g^{M}(\widetilde{e}_i,\widetilde{e}_j);~~\nonumber\\
&c(\widetilde{e}_i)\widehat{c}(\widetilde{e}_j)+\widehat{c}(\widetilde{e}_j)c(\widetilde{e}_i)=0.\nonumber\\
\end{align}
The Hodge-Dirac operator is given in \cite{Y}
\begin{align}
\label{a5}
\widetilde{D}=d+\delta=\sum^n_{i=1}c(\widetilde{e}_i)\bigg[\widetilde{e}_i+\frac{1}{4}\sum_{s,t}\omega_{s,t}
(\widetilde{e}_i)[\widehat{c}(\widetilde{e}_s)\widehat{c}(\widetilde{e}_t)
-c(\widetilde{e}_s)c(\widetilde{e}_t)]\bigg],
\end{align}
where $c(\widetilde{e}_i)$ and $\widehat{c}(\widetilde{e}_s)$ denote the Clifford action.

The following Lemma of Dabrowski etc.'s Einstein functional \cite{DL1} play a key role in our proof of the Einstein functional for manifold with boundary. Let $\overline{v}$ and $\overline{w}$ with the components with respect to local coordinates
$\overline{v}_a$ and $\overline{w}_b$, respectively, be two differential forms represented in such a way as endomorphisms (matrices)
$c(\overline{v})$ and $c(\overline{w})$ on $\Gamma(M,\bigwedge^*(T^*M))$. We assume thus that M is a $n=2m$ dimensional oriented Riemannian manifold and
use the Clifford action of one-forms as $0$-order differential operators. Using the operator $c(\overline{w})(\widetilde{D}c(\overline{v})+c(\overline{v})\widetilde{D})\widetilde{D}^{-n+1}$ acting on $\Gamma(M,\bigwedge^*(T^*M))$, the spectral functionals over the dual bimodule of one-forms defined by
\begin{lem}\label{lema1}\cite{DL1}The Einstein functional equals to
\begin{align}\label{a1}
{\rm Wres}[c(\overline{w})(\widetilde{D}c(\overline{v})+c(\overline{v})\widetilde{D})\widetilde{D}^{-n+1}]=\frac{v_{n-1}}{6}2^n\int_M[Ric^{ab}-\frac{1}{2}s(g)g^{ab}]\overline{v}_a\overline{w}_b{Vol_g},
\end{align}
where $g^*(\overline{v},\overline{w})=g^{ab}\overline{v}_a\overline{w}_b$ and $G(\overline{v},\overline{w})=(Ric^{ab}-\frac{1}{2}s(g)g^{ab})\overline{v}_a\overline{w}_b$ denote the Einstein tensor evaluated on two one-forms, where $\overline{v}=\sum_{j=1}^n\overline{v}_jdx_j,$ $\overline{w}=\sum_{l=1}^n\overline{w}_ldx_l,$ and $v_{n-1}=\frac{2\pi^{\frac{n}{2}}}{\Gamma(\frac{n}{2})}$.
\end{lem}
Let $\overline{v}=\sum_{j=1}^nv_je^{j,*}$, $\overline{w}=\sum_{l=1}^nw_le^{l,*}$, where $\{e^{1,*},e^{2,*},\cdot\cdot\cdot,e^{n,*}\},$ is the orthogonal basis about $g^{TM,*}$. Let $v=\sum_{j=1}^n\overline{v}(\widetilde{e}_j)\widetilde{e}_j:=\sum_{j=1}^nv_j\widetilde{e}_j,$ $w=\sum_{l=1}^nw_l\widetilde{e}_l,$ be the vector fields dual to one forms $\overline{v}$, $\overline{w}$. By the definition of $c(\overline{v})=c(v)$, $c(\overline{w})=c(w)$ and Lemma \ref{lema1}, we get
\begin{lem}\label{lema2}The Einstein functional equals to
\begin{align}\label{a2}
{\rm Wres}[c(w)(\widetilde{D}c(v)+c(v)\widetilde{D})\widetilde{D}^{-n+1}]=\frac{v_{n-1}}{6}2^n\int_M[Ric(v,w)-\frac{1}{2}s(g)g(v,w)]{Vol_g},
\end{align}
where $g(v,w)$ denotes the inner product evaluated on the two vector fields.
\end{lem}
Using an explicit formula for the spectral functionals of above the Hodge-Dirac operator, we can reformulate Theorems
for manifold $(M,g^M)$ with boundary $\partial M$ as follows
\begin{prop}\label{propa1}For the type-I operator, the Einstein functional for $4$ dimensional oriented Riemannian manifolds with
boundary equals to
\begin{align}
\label{a3}
\widetilde{{\rm Wres}}[\pi^+(c(w)(\widetilde{D}c(v)+c(v)\widetilde{D})\widetilde{D}^{-1})\circ\pi^+(\widetilde{D}^{-2})]={\rm Wres}[c(w)(\widetilde{D}c(v)+c(v)\widetilde{D})\widetilde{D}^{-3}]+\int_{\partial M}\Phi,
\end{align}
where
\begin{align}
\label{a4}
\Phi &=\int_{|\xi'|=1}\int^{+\infty}_{-\infty}\sum^{\infty}_{j, k=0}\sum\frac{(-i)^{|\alpha|+j+k+1}}{\alpha!(j+k+1)!}
{\rm tr}_{\wedge^*T^*M\bigotimes\mathbb{C}}[\partial^j_{x_n}\partial^\alpha_{\xi'}\partial^k_{\xi_n}\sigma^+_{r}(c(w)(\widetilde{D}c(v)+c(v)\widetilde{D})\widetilde{D}^{-1})(x',0,\xi',\xi_n)
\nonumber\\
&\times\partial^\alpha_{x'}\partial^{j+1}_{\xi_n}\partial^k_{x_n}\sigma_{l}(\widetilde{D}^{-2})(x',0,\xi',\xi_n)]d\xi_n\sigma(\xi')dx',
\end{align}
and the sum is taken over $r+l-k-j-|\alpha|=-3$.
\end{prop}
\begin{prop}\label{propa2}For the type-II operator, the Einstein functional for $4$ dimensional oriented Riemannian
 manifolds with
boundary equals to
\begin{align}
\label{a5}
\widetilde{{\rm Wres}}[\pi^+((c(w)(\widetilde{D}c(v)+c(v)\widetilde{D})\widetilde{D}^{-2})\circ\pi^+(\widetilde{D}^{-1})]={\rm Wres}[c(w)(\widetilde{D}c(v)+c(v)\widetilde{D})\widetilde{D}^{-3}]+\int_{\partial M}\Psi,
\end{align}
where
\begin{align}
\label{a6}
\Psi &=\int_{|\xi'|=1}\int^{+\infty}_{-\infty}\sum^{\infty}_{j, k=0}\sum\frac{(-i)^{|\alpha|+j+k+1}}{\alpha!(j+k+1)!}
 {\rm tr}_{\wedge^*T^*M\bigotimes\mathbb{C}}[\partial^j_{x_n}\partial^\alpha_{\xi'}\partial^k_{\xi_n}\sigma^+_{r}(c(w)(\widetilde{D}c(v)+c(v)\widetilde{D})\widetilde{D}^{-2})(x',0,\xi',\xi_n)
\nonumber\\
&\times\partial^\alpha_{x'}\partial^{j+1}_{\xi_n}\partial^k_{x_n}\sigma_{l}(\widetilde{D}^{-1})(x',0,\xi',\xi_n)]d\xi_n\sigma(\xi')dx',
\end{align}
and the sum is taken over $r+l-k-j-|\alpha|=-3$.
\end{prop}
\section{ The type-I operator and Dabrowski-Sitarz-Zalecki type theorems for 4-dimensional manifolds with boundary}
\label{section:4}
 In this section, we compute the type-I operator and prove the Dabrowski-Sitarz-Zalecki type theorems for 4-dimensional manifolds with boundary.
 
By Propsition 3.1 in \cite{WJ2}, we get the following propsition
\begin{prop}\label{propb1}For the type-I operator, the Einstein functional for $4$ dimensional oriented Riemannian manifolds with boundary is defined by
\begin{align}\label{b1}
&\widetilde{{\rm Wres}}[\pi^+(c(w)(\widetilde{D}c(v)+c(v)\widetilde{D})\widetilde{D}^{-1})\circ\pi^+(\widetilde{D}^{-2})]\nonumber\\
&=\widetilde{{\rm Wres}}[\pi^+(c(w)c(\widetilde{e}_j)c(\nabla^{T^*M}_{\widetilde{e}_j}v)\widetilde{D}^{-1})\circ\pi^+(\widetilde{D}^{-2})]+\widetilde{{\rm Wres}}[\pi^+(-2c(w)\nabla^{\bigwedge^*T^*M
}_v\widetilde{D}^{-1})\circ\pi^+(\widetilde{D}^{-2})].
\end{align}
\end{prop}
Since $\Phi^a$ and $\Phi^b$ are global forms on $\partial M$, so for any fixed point $x_0\in\partial M$, we choose the normal coordinates
$U$ of $x_0$ in $\partial M$ (not in $M$) and compute $\Psi(x_0)$ in the coordinates $\widetilde{U}=U\times [0,1)\subset M$ and the
metric $\frac{1}{h(x_n)}g^{\partial M}+dx_n^2.$ The dual metric of $g^M$ on $\widetilde{U}$ is ${h(x_n)}g^{\partial M}+dx_n^2.$  Write
$g^M_{ij}=g^M(\frac{\partial}{\partial x_i},\frac{\partial}{\partial x_j});~ g_M^{ij}=g^M(dx_i,dx_j)$, then

\begin{equation}
[g^M_{ij}]= \left[\begin{array}{lcr}
  \frac{1}{h(x_n)}[g_{ij}^{\partial M}]  & 0  \\
   0  &  1
\end{array}\right];~~~
[g_M^{ij}]= \left[\begin{array}{lcr}
  h(x_n)[g^{ij}_{\partial M}]  & 0  \\
   0  &  1
\end{array}\right],
\end{equation}
and
\begin{equation}
\partial_{x_s}g_{ij}^{\partial M}(x_0)=0, 1\leq i,j\leq n-1; ~~~g_{ij}^M(x_0)=\delta_{ij}.
\end{equation}
\indent Let $n=4$ and $\{e_1,\cdot\cdot\cdot,e_n\}$ be an orthonormal frame field in $U$ about $g^{\partial M}$ which is parallel along geodesics and $e_i=\frac{\partial}{\partial x_i}(x_0)$, then $\{\widetilde{e_1}=\sqrt{h(x_n)}e_1,\cdot\cdot\cdot,\widetilde{e_{n-1}}=\sqrt{h(x_n)}e_{n-1},\widetilde{e_n}=\frac{\partial}{\partial x_n}\}$ is the orthonormal frame field in $\widetilde{U}$ about $g^M$. Locally $S(TM)|\widetilde{U}\cong \widetilde{U}\times\wedge^*_C(\frac{n}{2})$. Let $\{f_1,\cdot\cdot\cdot,f_n\}$ be the orthonormal basis of $\wedge^*_C(\frac{n}{2})$. Take a spin frame field $\sigma:\widetilde{U}\rightarrow Spin(M)$ such that $\pi\sigma=\{\widetilde{e_1},\cdot\cdot\cdot,\widetilde{e}_n\}$ where $\pi:Spin(M)\rightarrow O(M)$ is a double covering, then $\{[\sigma,f_i],1\leq i\leq4\}$ is an orthonormal frame of $S(TM)|_{\widetilde{U}}$. In the following, since the global form $\Psi$ is independent of the choice of the local frame, so we can compute $tr_{S(TM)}$ in the frame $\{[\sigma,f_i],1\leq i\leq4\}$. Let $\{\widehat{e_1},\cdot\cdot\cdot,\widehat{e}_n\}$ be the canonical basis of $R^n$ and $c(\widehat{e_i})\in H om(\wedge^*_C(\frac{n-1}{2}),\wedge^*_C(\frac{n-1}{2}))$ be the Clifford action, then
\begin{align}
c(\widetilde{e_i})=[\sigma,c(\widehat{e_i})];~~~c(\widetilde{e_i})[\sigma,f_i]=[\sigma,c(\widehat{e_i})f_i];~~~\frac{\partial}{\partial x_i}=[\sigma,\frac{\partial}{\partial x_i}].
\end{align}
then we have $\frac{\partial}{\partial x_i}c(\widetilde{e_i})=0$ in the above frame.
\subsection{$\widetilde{{\rm Wres}}[\pi^+(c(w)c(\widetilde{e}_j)c(\nabla^{T*M}_{\widetilde{e}_j}v)\widetilde{D}^{-1})\circ\pi^+(\widetilde{D}^{-2})]$}
Combining with the generating Proposition \ref{propb1}, this yields
\begin{align}\label{b2}
\widetilde{{\rm Wres}}[\pi^+(c(w)c(\widetilde{e}_j)c(\nabla^{T^*M}_{\widetilde{e}_j}v)\widetilde{D}^{-1})\circ\pi^+(\widetilde{D}^{-2})]={\rm Wres}[(c(w)c(\widetilde{e}_j)c(\nabla^{T^*M}_{\widetilde{e}_j}v)\widetilde{D}^{-3}]+\int_{\partial M}\Phi^a.\nonumber\\
\end{align}

Let $\xi=\sum_{j}\xi_{j}dx_{j}$ and $\nabla^L_{\partial_{i}}\partial_{j}=\sum_{k}\Gamma_{ij}^{k}\partial_{k}$,  we denote that
\begin{align}
\label{a19}
&\sigma_{i}=-\frac{1}{4}\sum_{s,t}\omega_{s,t}
(\widetilde{e}_i)c(\widetilde{e}_s)c(\widetilde{e}_t)
;~~~a_{i}=\frac{1}{4}\sum_{s,t}\omega_{s,t}
(\widetilde{e}_i)\widehat{c}(\widetilde{e}_s)\widehat{c}(\widetilde{e}_t);\nonumber\\
&\xi^{j}=g^{ij}\xi_{i};~~~~\Gamma^{k}=g^{ij}\Gamma_{ij}^{k};~~~~\sigma^{j}=g^{ij}\sigma_{i};
~~~~a^{j}=g^{ij}a_{i}.
\end{align}
Then, we get the following Lemma
\begin{lem}\label{lemb1} The following identities hold:
\begin{align}\label{b3}
&\sigma_{-1}(\widetilde{D}^{-1})=\frac{ic(\xi)}{|\xi|^2};\nonumber\\
&\sigma_{-2}(\widetilde{D}^{-2})=|\xi|^{-2};\nonumber\\
&\sigma_{-2}(\widetilde{D}^{-1})=\frac{c(\xi)\sigma_{0}(\widetilde{D})c(\xi)}{|\xi|^4}+\frac{c(\xi)}{|\xi|^6}\sum_jc(dx_j)
\Big[\partial_{x_j}(c(\xi))|\xi|^2-c(\xi)\partial_{x_j}(|\xi|^2)\Big];\nonumber\\
&\sigma_{-3}(\widetilde{D}^{-2})=-\sqrt{-1}|\xi|^{-4}\xi_k(\Gamma^k-2\sigma^k+2a^k)-\sqrt{-1}|\xi|^{-6}2\xi^j\xi_\alpha\xi_\beta\partial_jg^{\alpha\beta}.\nonumber\\
\end{align}
\end{lem}
Now, we can compute $\Phi^a$. When $n=4$, then ${\rm tr}_{\bigwedge^*T^*M} [{\rm \texttt{id}}]=2^4$, the sum is taken over $
r+l-k-j-|\alpha|=-3,~~r=-1,~~l=-2,$ then we have the only one case:
\begin{align}\label{b4}
\Phi^a=-i\int_{|\xi'|=1}\int^{+\infty}_{-\infty}
 {\rm tr}[\pi^+_{\xi_n}\sigma_{-1}(\sum_{j=1}^nc(w)c(\widetilde{e}_j)c(\nabla^{T*M}_{\widetilde{e}_j}v)\widetilde{D}^{-1})\times\partial_{\xi_n}\sigma_{-2}(\widetilde{D}^{-2})](x_0)d\xi_n\sigma(\xi')dx'.
\end{align}
By Lemma \ref{lemb1}, we get
\begin{align}
\pi^+_{\xi_n}\sigma_{-1}\left(\sum_{j=1}^nc(w)c(\widetilde{e}_j)c(\nabla^{T^*M}_{\widetilde{e}_j}v)\widetilde{D}^{-1}\right)=\sum_{j=1}^nc(w)c(\widetilde{e}_j)c(\nabla^{T^*M}_{\widetilde{e}_j}v)\frac{c(\xi')+ic(dx_n)}{2(\xi_n-i)},\nonumber
\end{align}
and by further derivative, we get
\begin{align} 
\partial_{\xi_n}\sigma_{-2}(\widetilde{D}^{-2})=\frac{-2\xi_n}{(1+\xi_n^2)^2}.\nonumber
\end{align}
It follows that 
\begin{align}
&{\rm tr}[\pi^+_{\xi_n}\sigma_{-1}(\sum_{j=1}^nc(w)c(\widetilde{e}_j)c(\nabla^{T^*M}_{\widetilde{e}_j}v)\widetilde{D}^{-1})\times\partial_{\xi_n}\sigma_{-2}(\widetilde{D}^{-2})](x_0)\nonumber\\
&=-\frac{\xi_n}{(\xi_n-i)(1+\xi_n^2)^2}{\rm tr}\sum_{j=1}^n[c(w)c(\widetilde{e}_j)c(\nabla^{T^*M}_{\widetilde{e}_j}v)c(\xi')](x_0)\nonumber\\
&-\frac{i\xi_n}{(\xi_n-i)(1+\xi_n^2)^2}{\rm tr}\sum_{j=1}^n[c(w)c(\widetilde{e}_j)c(\nabla^{T^*M}_{\widetilde{e}_j}v)c(dx_n)](x_0)\nonumber.
\end{align}
By the relation of the Clifford action and ${\rm tr}{ab}={\rm tr }{ba}$,  we have the equalities:
\begin{align}\label{123}
{\rm tr}\sum_{j=1}^n[c(w)c(\widetilde{e}_j)c(\nabla^{T*M}_{\widetilde{e}_j}v)c(dx_n)](x_0)=\bigg(\sum_{j=1}^ng(e_j,\nabla^L_{e_j}v)g(w,\frac{\partial}{\partial{x_n}})-g(w,\nabla^L_{\frac{\partial}{\partial{x_n}}}v)+g(\nabla^L_wv,\frac{\partial}{\partial{x_n}})\bigg) {\rm {\rm tr}[\texttt{id}]}.
\end{align}
We note that $i<n,~\int_{|\xi'|=1}\xi_{i_{1}}\xi_{i_{2}}\cdots\xi_{i_{2d+1}}\sigma(\xi')=0$,
so ${\rm tr}\sum_{j=1}^n[c(w)c(\widetilde{e}_j)c(\nabla^{T*M}_{\widetilde{e}_j}v)c(\xi')](x_0)$ has no contribution for computing $\Phi^a$. Then, we have

\begin{align}\label{iiii}
\Phi^a&=-i\int_{|\xi'|=1}\int^{+\infty}_{-\infty}
-\frac{i\xi_n}{(\xi_n-i)(1+\xi_n^2)^2}{\rm tr}\sum_{j=1}^n[c(w)c(\widetilde{e}_j)c(\nabla^{T^*M}_{\widetilde{e}_j}v)c(dx_n)](x_0)d\xi_n\sigma(\xi')dx'\nonumber\\
&=-\int_{\Gamma^+}\frac{\xi_n}{(\xi_n-i)(1+\xi_n^2)^2}d\xi_n\bigg(\sum_{j=1}^ng(e_j,\nabla^L_{e_j}v)g(w,\frac{\partial}{\partial{x_n}})-g(w,\nabla^L_{\frac{\partial}{\partial{x_n}}}v)+g(\nabla^L_wv,\frac{\partial}{\partial{x_n}})\bigg){\rm {\rm tr}[\texttt{id}]}\Omega_3dx'\nonumber\\
&=-\frac{2\pi i}{2!}\left[\frac{\xi_n}{(\xi_n+i)^2}\right]^{(2)}\bigg|_{\xi_n=i}\bigg(\sum_{j=1}^ng(e_j,\nabla^L_{e_j}v)g(w,\frac{\partial}{\partial{x_n}})-g(w,\nabla^L_{\frac{\partial}{\partial{x_n}}}v)+g(\nabla^L_wv,\frac{\partial}{\partial{x_n}})\bigg)\times16\times4\pi dx'\nonumber\\
&=-8\bigg(\sum_{j=1}^ng(e_j,\nabla^L_{e_j}v)g(w,\frac{\partial}{\partial{x_n}})-g(w,\nabla^L_{\frac{\partial}{\partial{x_n}}}v)+g(\nabla^L_wv,\frac{\partial}{\partial{x_n}})\bigg)\pi^2dx'.
\end{align}

\subsection{$\widetilde{{\rm Wres}}[\pi^+(-2c(w)\nabla^{\bigwedge^*T^*M} _v\widetilde{D}^{-1})\circ\pi^+(\widetilde{D}^{-2})]$}
Combining with the generating Proposition \ref{propb1}, this yields
\begin{align}
\widetilde{{\rm Wres}}[\pi^+(-2c(w)\nabla^{\bigwedge^*T^*M} _v\widetilde{D}^{-1})\circ\pi^+(\widetilde{D}^{-2})]={\rm Wres}[-2c(w)\nabla^{\bigwedge^*T^*M} _v\widetilde{D}^{-3}]+\int_{\partial M}\Phi^b.\nonumber
\end{align}
 We define $\nabla_v^{\bigwedge^*T^*M} :=v+\frac{1}{4}\sum_{ij}\langle\nabla_v^L{\widetilde{e}_i},\widetilde{e}_j\rangle [c(\widetilde{e}_i)c(\widetilde{e}_j)-\widehat{c}(\widetilde{e}_i)\widehat{c}(\widetilde{e}_j)]$. Set $A(v)=\frac{1}{4}\sum_{ij}\langle\nabla_v^L{\widetilde{e}_i},\widetilde{e}_j\rangle [c(\widetilde{e}_i)c(\widetilde{e}_j)-\widehat{c}(\widetilde{e}_i)\widehat{c}(\widetilde{e}_j)]$, where $v=\sum_{j=1}^nv_j\partial_{x_j}$ and $\sigma(\partial_{x_j})=\sqrt{-1}\xi_j$. Let $v=v^T+v_n\partial_n,~w=w^T+w_n\partial_n,$ then we have $\sum_{j=1}^{n-1}v_jw_j=g(v^T,w^T)(x_0).$ Therefore, we have the following lemma
\begin{lem}\label{lem553} The following identities hold:
\begin{align}
\label{b42}
&\sigma_{0}(\nabla_v^{\bigwedge^*T^*M})=A(v);\nonumber\\
&\sigma_{1}(\nabla_v^{\bigwedge^*T^*M})=\sqrt{-1}\sum_{j=1}^nv_j\xi_j.\nonumber
\end{align}
\end{lem}
\indent Write
 \begin{eqnarray}
\widetilde{D}_x^{\alpha}&=(-i)^{|\alpha|}\partial_x^{\alpha};
~\sigma(\widetilde{D})=p_1+p_0;
~\sigma(\widetilde{D}^{-1})=\sum^{\infty}_{j=1}q_{-j}.
\end{eqnarray}

\indent By the composition formula of pseudodifferential operators, we have
\begin{align*}
1=\sigma(\widetilde{D}\circ \widetilde{D}^{-1})&=\sum_{\alpha}\frac{1}{\alpha!}\partial^{\alpha}_{\xi}[\sigma(\widetilde{D})]
\widetilde{D}_x^{\alpha}[\sigma(\widetilde{D}^{-1})]\nonumber\\
&=(p_1+p_0)(q_{-1}+q_{-2}+q_{-3}+\cdots)\nonumber\\
&~~~+\sum_j(\partial_{\xi_j}p_1+\partial_{\xi_j}p_0)(
\widetilde{D}_{x_j}q_{-1}+\widetilde{D}_{x_j}q_{-2}+\widetilde{D}_{x_j}q_{-3}+\cdots)\nonumber\\
&=p_1q_{-1}+(p_1q_{-2}+p_0q_{-1}+\sum_j\partial_{\xi_j}p_1\widetilde{D}_{x_j}q_{-1})+\cdots,
\end{align*}
so
\begin{equation*}
q_{-1}=p_1^{-1};~q_{-2}=-p_1^{-1}[p_0p_1^{-1}+\sum_j\partial_{\xi_j}p_1\widetilde{D}_{x_j}(p_1^{-1})].
\end{equation*}
 \begin{lem} \cite{Wa3}\label{le:32}
With the metric $g^{M}$ on $M$ near the boundary
\begin{align}
\label{b1q1}
\partial_{x_j}(|\xi|_{g^M}^2)(x_0)&=\left\{
       \begin{array}{c}
        0,  ~~~~~~~~~~ ~~~~~~~~~~ ~~~~~~~~~~~~~{\rm if }~j<n, \\[2pt]
       h'(0)|\xi'|^{2}_{g^{\partial M}},~~~~~~~~~~~~~~~~~~~~{\rm if }~j=n,
       \end{array}
    \right. \nonumber\\
\partial_{x_j}[c(\xi)](x_0)&=\left\{
       \begin{array}{c}
      0,  ~~~~~~~~~~ ~~~~~~~~~~ ~~~~~~~~~~~~~{\rm if }~j<n,\\[2pt]
\frac{1}{2}h'(0)c(\xi')(x_{0}), ~~~~~~~~~~~~~~~~~{\rm if }~j=n,
       \end{array}
    \right.
\end{align}
where $\xi=\xi'+\xi_{n}dx_{n}$.
\end{lem}
By Lemma \ref{lem553} and Lemma \ref{le:32}, we have the following lemma
\begin{lem}\label{lemb2e}The following identities hold:
\begin{align}
\sigma_{0}(\nabla_v^{\bigwedge^*T^*M} \widetilde{D}^{-1})&=-\sum_{j=1}^nv_j\xi_j\frac{c(\xi)}{|\xi|^2};\nonumber\\
\sigma_{-1}(\nabla_v^{\bigwedge^*T^*M} \widetilde{D}^{-1})&=\sigma_1(\nabla_v^{\bigwedge^*T^*M} )\sigma_{-2}(\widetilde{D}^{-1})+\sigma_{0}(\nabla_v^{\bigwedge^*T^*M} )\sigma_{-1}(\widetilde{D}^{-1})+\sum_{j=1}^n\xi_j\sigma_{1}(\nabla_v^{\bigwedge^*T^*M} )D_{x_j}\sigma_{-1}(\widetilde{D}^{-1})\nonumber\\
&=\sqrt{-1}\sum_{j=1}^nv_j\xi_j\bigg[\frac{c(\xi)\sigma_{0}(\widetilde{D})c(\xi)}{|\xi|^4}+\frac{c(\xi)}{|\xi|^6}\sum_jc(dx_j)
\Big(\partial_{x_j}(c(\xi)|\xi|^2-c(\xi)\partial_{x_j}(|\xi|^2)\Big)\bigg]\nonumber\\
&+A(v)\frac{\sqrt{-1}c(\xi)}{|\xi|^2}+v_n\bigg(\frac{\sqrt{-1}\partial_{x_n}c(\xi')(x_0)}{|\xi|^2}-\frac{\sqrt{-1}c(\xi)|\xi'|^2h'(0)}{|\xi|^4}\bigg).\nonumber
\end{align}
\end{lem}
From Lemma \ref{lema2}, we obtain the following result
\begin{align}\label{llll}
{\rm Wres}[c(w)c(\widetilde{e}_j)c(\nabla^{T^*M}_{\widetilde{e}_j}v)\widetilde{D}^{-3}]+{\rm Wres}[-2c(w)\nabla^{\bigwedge^*T^*M} _v\widetilde{D}^{-3}]&={\rm Wres}[c(w)(\widetilde{D}c(v)+c(v)\widetilde{D})\widetilde{D}^{-3}]\nonumber\\
&=\frac{64\pi^2}{3}\int_M[Ric(v,w)-\frac{1}{2}s(g)g(v,w)]{Vol_g}.
\end{align}
 
Therefore, we only need to compute $\Phi^b$. The sum is taken over $
r+l-k-j-|\alpha|=-3,~~r\leq 0,~~l\leq -2,$ then we have the following five cases:
~\\
\noindent  {\bf case a)~I)}~$r=0,~l=-2,~k=j=0,~|\alpha|=1$.\\
\noindent By (\ref{a3}), we get
\begin{align}
\Phi^b_1=-\int_{|\xi'|=1}\int^{+\infty}_{-\infty}\sum_{|\alpha|=1}
 {\rm tr}[\partial^\alpha_{\xi'}\pi^+_{\xi_n}\sigma_{0}(-2c(w)\nabla^{S(TM)}_v\widetilde{D}^{-1})\times
 \partial^\alpha_{x'}\partial_{\xi_n}\sigma_{-2}(\widetilde{D}^{-2})](x_0)d\xi_n\sigma(\xi')dx'.\nonumber
\end{align}
By Lemma 2.2 in \cite{Wa3}, for $i<n$, then
\begin{align}
\partial_{x_i}\sigma_{-2}(\widetilde{D}^{-2})(x_0)=\partial_{x_i}{(|\xi|^{-2})}(x_0)
=\frac{\partial_{x_i}(|\xi|^2)(x_0)}{|\xi|^4}=0,\nonumber
\end{align}
so $\Phi^b_1=0$.\\
 \noindent  {\bf case a)~II)}~$r=0,~l=-2,~k=|\alpha|=0,~j=1$.\\
\noindent By (\ref{a3}), we get
\begin{align}
\Phi^b_2=-\frac{1}{2}\int_{|\xi'|=1}\int^{+\infty}_{-\infty} {\rm
tr} [\partial_{x_n}\pi^+_{\xi_n}\sigma_{0}(-2c(w)\nabla^{\bigwedge^*T^*M} _v\widetilde{D}^{-1})\times
\partial_{\xi_n}^2\sigma_{-2}(\widetilde{D}^{-2})](x_0)d\xi_n\sigma(\xi')dx'.
\end{align}
\noindent By Lemma \ref{lemb1}, we have\\
\begin{align}
\partial_{\xi_n}^2\sigma_{-2}(\widetilde{D}^{-2})(x_0)=\partial_{\xi_n}^2(|\xi|^{-2})(x_0)=\frac{6\xi_n^2-2}{(1+\xi_n^2)^3}.\nonumber
\end{align}
By Lemma \ref{lemb2e}, we have\\
\begin{align}\label{b27}
\partial_{x_n}\sigma_{0}(-2c(w)\nabla_v^{S(TM)}\widetilde{D}^{-1})&=2\sum_{j=1}^{n-1}\xi_j\bigg[\frac{1}{|\xi|^2}\partial_{x_n}(v_jc(w)c(\xi'))-\frac{h'(0)|\xi'|^2}{|\xi|^4}v_jc(w)c(\xi')+\frac{\xi_n}{|\xi|^2}\partial_{x_n}(v_jc(w)c(\xi'))\nonumber\\
&-\frac{h'(0)|\xi'|^2}{|\xi|^4}v_jc(w)c(dx_n)\bigg]+\frac{2\xi_n}{|\xi|^2}\partial_{x_n}(v_nc(w)c(\xi'))-\frac{2\xi_nh'(0)|\xi'|^2}{|\xi|^4}v_nc(w)c(\xi')\nonumber\\
&+\frac{2\xi_n^2}{|\xi|^2}\partial_{x_n}(v_nc(w)c(\xi'))-\frac{2\xi_nh'(0)|\xi'|^2}{|\xi|^4}v_nc(w)c(dx_n).
\end{align}
We note that $i<n,~\int_{|\xi'|=1}\xi_{i_{1}}\xi_{i_{2}}\cdots\xi_{i_{2d+1}}\sigma(\xi')=0$,
so we omit some items that have no contribution for computing $\Phi^b_2$. Moreover
\begin{align}\label{b28}
&\partial_{x_n}\pi^+_{\xi_n}\sigma_{0}(-2c(w)\nabla_v^{\bigwedge^*T^*M} \widetilde{D}^{-1})\nonumber\\
&=\pi^+_{\xi_n}\partial_{x_n}\sigma_{0}(-2c(w)\nabla_v^{\bigwedge^*T^*M} \widetilde{D}^{-1})\nonumber\\
&=\sum_{j=1}^{n-1}\xi_j\bigg[\frac{-i}{(\xi_n-i)}\partial_{x_n}(v_jc(w)c(\xi'))+\frac{2+i\xi_n}{2(\xi_n-i)^2}h'(0)|\xi'|^2v_jc(w)c(\xi')\bigg]+\frac{2\xi_n}{|\xi|^2}\partial_{x_n}(v_jc(w)c(\xi'))\nonumber\\
&-\frac{2\xi_nh'(0)|\xi'|^2}{|\xi|^4}v_jc(w)c(\xi')+\frac{i}{\xi_n-i}\partial_{x_n}(v_nc(w)c(dx_n))+\frac{i\xi_nh'(0)|\xi'|^2}{2(\xi_n-i)^2}h'(0)|\xi'|^2v_nc(w)c(dx_n).
\end{align}
By the relation of the Clifford action and ${\rm tr}{ab}={\rm tr }{ba}$,  we have the equalities:
\begin{align}\label{eee}
&\sum_{j=1}^{n-1}\xi_j{\rm tr }[\partial_{x_n}(v_jc(w)c(\xi'))]=\bigg(-\sum_{j,k=1}^{n-1}\xi_j\xi_k\partial_{x_n}(v_jw_k)-\frac{1}{2}\sum_{j,k=1}^{n-1}\xi_j\xi_kh'(0)v_jw_k\bigg){\rm tr}[\texttt{id}];\nonumber\\
&\sum_{j=1}^{n-1}\xi_jv_j{\rm tr }[c(w)c(\xi')]=-\sum_{j,k=1}^{n-1}\xi_j\xi_kv_jw_k{\rm tr}[\texttt{id}];~~~~~{\rm tr }[\partial_{x_n}(v_nc(w)c(dx_n))]=-\partial_{x_n}(v_nw_n){\rm tr}[\texttt{id}];\nonumber\\
&{\rm tr }[v_nc(w)c(dx_n)]=-v_nw_n{\rm tr}[\texttt{id}].
\end{align}
Then, we have
\begin{align}\label{b78}
& {\rm
tr} [\partial_{x_n}\pi^+_{\xi_n}\sigma_{0}(-2c(w)\nabla^{\bigwedge^*T^*M} _v\widetilde{D}^{-1})\times
\partial_{\xi_n}^2\sigma_{-2}(\widetilde{D}^{-2})](x_0)\nonumber\\
&=\frac{i(2-6\xi_n^2)}{(\xi_n-i)(1+\xi_n^2)^3}\bigg(-\sum_{j,k=1}^{n-1}\xi_j\xi_k\partial_{x_n}(v_jw_k)-\frac{1}{2}\sum_{j,k=1}^{n-1}\xi_j\xi_kh'(0)v_jw_k\bigg){\rm tr}[\texttt{id}]\nonumber\\
&-\frac{(2+i\xi_n)(3\xi_n^2-1)}{(\xi_n-i)^2(1+\xi_n^2)^3}h'(0)|\xi'|^2\sum_{j,k=1}^{n-1}\xi_j\xi_kv_jw_k{\rm tr}[\texttt{id}]-\frac{i(6\xi_n^2-2)}{(\xi_n-i)(1+\xi_n^2)^3}\partial_{x_n}(v_nw_n){\rm tr}[\texttt{id}]\nonumber\\
&-\frac{i\xi_n(3\xi_n^2-1)}{(\xi_n-i)^2(1+\xi_n^2)^3}h'(0)|\xi'|^2v_nw_n{\rm tr}[\texttt{id}].
\end{align}
Therefore, by $\int_{|\xi'|=1}\xi_j\xi_k=\frac{4\pi}{3}\delta_j^k,$ we have
\begin{align}\label{35}
\Phi^b_2&=-\frac{1}{2}\int_{|\xi'|}\int_{-\infty}^{\infty}{\rm
tr} [\partial_{x_n}\pi^+_{\xi_n}\sigma_{0}(-2c(w)\nabla^{\bigwedge^*T^*M} _v\widetilde{D}^{-1})\times
\partial_{\xi_n}^2\sigma_{-2}(\widetilde{D}^{-2})](x_0)d\xi_n\sigma(\xi')dx'\nonumber\\
&=-\frac{1}{2}\times\frac{2\pi i}{3!}\left[\frac{i(2-6\xi_n)}{(\xi_n+i)^3}\right]^{(3)}\bigg|_{\xi_n=i}\times -\frac{4\pi}{3}\sum_{j=1}^{n-1}[\partial_{x_n}(v_jw_j)+\frac{1}{2}h'(0)v_jw_j]\times 16dx'\nonumber\\
&-\frac{1}{2}\times\frac{2\pi i}{4!}\left[\frac{3i\xi_n^3+6\xi_n^2-i\xi_n-2}{(\xi_n+i)^3}\right]^{(4)}\bigg|_{\xi_n=i}\times -\frac{4\pi}{3}\sum_{j=1}^{n-1}v_jw_j\times 16dx'\nonumber\\
&-\frac{1}{2}\times\frac{2\pi i}{3!}\left[\frac{i(6\xi_n^2-2)}{(\xi_n+i)^3}\right]^{(3)}\bigg|_{\xi_n=i}\times -\partial_{x_n}(v_nw_n)\times 16\Omega_3dx'\nonumber\\
&-\frac{1}{2}\times\frac{2\pi i}{4!}\left[\frac{3i\xi_n^4-i\xi_n^2}{(\xi_n+i)^3}\right]^{(4)}\bigg|_{\xi_n=i}\times -v_nw_n\times 16\Omega_3dx'\nonumber\\
 &=\bigg(-\frac{8}{3}\partial_{x_n}g(v^T,w^T)+2h'(0)g(v^T,w^T)+8\partial_{x_n}(v_nw_n)-2h'(0)v_nw_n\bigg)\pi^2dx'.
 \end{align}
\noindent  {\bf case a)~III)}~$r=0,~l=-2,~j=|\alpha|=0,~k=1$.\\
\noindent By (\ref{a3}), we get
\begin{align}\label{36}
\Phi^b_3&=-\frac{1}{2}\int_{|\xi'|=1}\int^{+\infty}_{-\infty}
{\rm tr} [\partial_{\xi_n}\pi^+_{\xi_n}\sigma_{0}(-2c(w)\nabla^{\bigwedge^*T^*M} _v\widetilde{D}^{-1})\times
\partial_{\xi_n}\partial_{x_n}\sigma_{-2}(\widetilde{D}^{-2})](x_0)d\xi_n\sigma(\xi')dx'.
\end{align}
By Lemma \ref{lemb1}, we have 
\begin{align}\label{37}
\partial_{x_n}\sigma_{-2}(D^{-2})(x_0)|_{|\xi'|=1}=\frac{4\xi_nh'(0)}{(1+\xi_n^2)^3}.
\end{align}
And by further calculation, we have
\begin{align}\label{38887}
\partial_{\xi_n}\partial_{x_n}\sigma_{-2}(D^{-2})(x_0)|_{|\xi'|=1}=-\frac{h'(0)}{(1+\xi_n^2)^2}.
\end{align}
By the Cauchy integral formula, we obtain
\begin{eqnarray}\label{64}
\pi^+_{\xi_n} \Big( \frac{1}{1+\xi_n^2}\Big)
&=&\frac{1}{2\pi i} \int_{\Gamma^+}\frac{1}{(\xi_n-\eta_n)(1+\xi_n^2)}d\eta_n\nonumber\\
&=&\left[\frac{1}{(\xi_n-\eta_n)(\eta_n+i)}\right]\Big|_{\eta_n=i}
     = \frac{-i}{2(\xi_n-i)},\nonumber\\
\pi^+_{\xi_n} \Big( \frac{\xi_n^2}{1+\xi_n^2}\Big)
&=&\frac{1}{2\pi i} \int_{\Gamma^+}\frac{\eta_n^2}{(\xi_n-\eta_n)(1+\xi_n^2)}d\eta_n\nonumber\\
&=&\frac{1}{2} \left[\frac{\eta_n}{(\xi_n-\eta_n)(\eta_n+i)}\right]\Big|_{\eta_n=i}
     = \frac{i}{2(\xi_n-i)}.
\end{eqnarray}
Then, we get
\begin{align}\label{mmmmm}
\partial_{\xi_n}\pi^+_{\xi_n}\sigma_{0}(-2c(w)\nabla^{\bigwedge^*T^*M} _v\widetilde{D}^{-1})
&=\sum_{j=1}^{n-1}v_j\xi_j\frac{i}{(\xi_n-i)^2}c(w)c(\xi')-v_n\frac{i}{(\xi_n-i)^2}c(w)c(dx_n).
\end{align}
By (\ref{eee}), we have
\begin{align}\label{kkk}
&{\rm tr} [\partial_{\xi_n}\pi^+_{\xi_n}\sigma_{0}(-2c(w)\nabla^{\bigwedge^*T^*M} _v\widetilde{D}^{-1})\times
\partial_{\xi_n}\partial_{x_n}\sigma_{-2}(\widetilde{D}^{-2})](x_0)\nonumber\\
&=-\frac{4i\xi_n}{(\xi_n-i)^2(1+\xi_n^2)^3}h'(0)\sum_{j,k=1}^{n-1}v_jw_k\xi_j\xi_k{\rm tr}[\texttt{id}]+\frac{4i\xi_n}{(\xi_n-i)^2(1+\xi_n^2)^3}h'(0)v_nw_n{\rm tr}[\texttt{id}].
\end{align}
Next, we perform the corresponding integral calculation on the above results. Therefore
\begin{align}\label{3ll5}
\Phi^b_3&=-\frac{1}{2}\int_{|\xi'|=1}\int^{+\infty}_{-\infty}
{\rm tr} [\partial_{\xi_n}\pi^+_{\xi_n}\sigma_{0}(-2c(w)\nabla^{\bigwedge^*T^*M} _v\widetilde{D}^{-1})\times
\partial_{\xi_n}\partial_{x_n}\sigma_{-2}(\widetilde{D}^{-2})](x_0)d\xi_n\sigma(\xi')dx'\nonumber\\
&=-\frac{1}{2}\times\frac{2\pi i}{4!}\left[\frac{4i\xi_n}{(\xi_n+i)^3}\right]^{(4)}\bigg|_{\xi_n=i}\times-\frac{4\pi}{3}\sum_{j=1}^{n-1}v_jw_j\times 16h'(0)dx'\nonumber\\
&-\frac{1}{2}\times\frac{2\pi i}{4!}\left[\frac{4i\xi_n}{(\xi_n+i)^3}\right]^{(4)}\bigg|_{\xi_n=i}\times -v_nw_n\times 16h'(0)\Omega_3dx'\nonumber\\
 &=\bigg(-\frac{10}{3}g(v^T,w^T)+10v_nw_n\bigg)h'(0)\pi^2dx'.
 \end{align}
\noindent  {\bf case b)}~$r=0,~l=-3,~k=j=|\alpha|=0$.\\
\noindent By (\ref{a3}), we get
\begin{align}\label{42}
\Phi^b_4&=-i\int_{|\xi'|=1}\int^{+\infty}_{-\infty}{\rm tr} [\pi^+_{\xi_n}\sigma_{0}(-2c(w)\nabla^{\bigwedge^*T^*M} _v\widetilde{D}^{-1})\times
\partial_{\xi_n}\sigma_{-3}(\widetilde{D}^{-2})](x_0)d\xi_n\sigma(\xi')dx'\nonumber\\
&=i\int_{|\xi'|=1}\int^{+\infty}_{-\infty}{\rm tr} [\partial_{\xi_n}\pi^+_{\xi_n}\sigma_{0}(-2c(w)\nabla^{\bigwedge^*T^*M} _v\widetilde{D}^{-1})\times
\sigma_{-3}(\widetilde{D}^{-2})](x_0)d\xi_n\sigma(\xi')dx'.
\end{align}

In the normal coordinate, $g^{ij}(x_{0})=\delta^{j}_{i}$ and $\partial_{x_{j}}(g^{\alpha\beta})(x_{0})=0$, if $j<n$; $\partial_{x_{j}}(g^{\alpha\beta})(x_{0})=h'(0)\delta^{\alpha}_{\beta}$, if $j=n$.
So by \cite{Wa3}, when $k<n$, we have $\Gamma^{n}(x_{0})=\frac{5}{2}h'(0)$, $\Gamma^{k}(x_{0})=0$, $\sigma^{n}(x_{0})=0$ and $\sigma^{k}=\frac{1}{4}h'(0)c(\widetilde{e}_{k})c(\widetilde{e}_{n}), a^{k}=-\frac{1}{4}h'(0)\widehat{c}(\widetilde{e}_{k})\widehat{c}(\widetilde{e}_{n})$. Then by Lemma \ref{lemb1}, we obtain
\begin{align}\label{43}
&\sigma_{-3}(D^{-2})(x_0)|_{|\xi'|=1}\nonumber\\
&=\frac{i}{(1+\xi_n^2)^2}
   \Big(\frac{1}{2}h'(0)\sum_{k<n}\xi_kc(\widetilde{e}_k)c(\widetilde{e}_n)+\frac{1}{2}h'(0)\sum_{k<n}\xi_k\widehat{c}(\widetilde{e}_k)\widehat{c}(\widetilde{e}_n)+\frac{5}{2}h'(0)\xi_n \Big)-\frac{2ih'(0)\xi_n}{(1+\xi_n^2)^3}\nonumber\\
   &=\frac{i}{2(1+\xi_n^2)^2}
   h'(0)\sum_{k<n}\xi_kc(\widetilde{e}_k)c(\widetilde{e}_n)+\frac{i}{2(1+\xi_n^2)^2}h'(0)\sum_{k<n}\xi_k\widehat{c}(\widetilde{e}_k)\widehat{c}(\widetilde{e}_n)-\frac{5i\xi_n^3+9i\xi_n}{2(1+\xi_n^2)^3}h'(0).
\end{align}
 By Lemma \ref{lemb2e}, we have
\begin{align}\label{45}
\partial_{\xi_n}\pi^+_{\xi_n}\sigma_{0}(-2c(w)\nabla_v^{\bigwedge^*T^*M} \widetilde{D}^{-1})&=\sum_{j=1}^{n-1}v_j\xi_j\bigg(\frac{i}{(\xi_n-i)^2}c(w)c(\xi')-\frac{1}{\xi_n-i)^2}c(w)c(dx_n)\bigg)\nonumber\\
&-v_n\bigg(\frac{1}{(\xi_n-i)^2}c(w)c(\xi')+\frac{i}{\xi_n-i)^2}c(w)c(dx_n)\bigg).
\end{align}
We note that $i<n,~\int_{|\xi'|=1}\xi_{i_{1}}\xi_{i_{2}}\cdots\xi_{i_{2d+1}}\sigma(\xi')=0$,
so we omit some items that have no contribution for computing $\Phi^b_2$.
By the relation of the Clifford action, we have the following identities:
\begin{align}\label{lll32}
&\sum_{k<n}\xi_k{\rm tr }[c(w)c(\widetilde{e}_k)c(\widetilde{e}_n)c(dx_n)]=\sum_{k=1}^{n-1}\xi_kw_k{\rm tr}[\texttt{id}];~~\sum_{k<n}\xi_k{\rm tr }[c(w)c(\widetilde{e}_k)c(\widetilde{e}_n)c(\xi')]=-\sum_{k,j=1}^{n-1}\xi_k\xi_jw_n{\rm tr}[\texttt{id}];\nonumber\\
&\sum_{k<n}\xi_k{\rm tr }[c(w)\widehat{c}(\widetilde{e}_k)\widehat{c}(\widetilde{e}_n)c(dx_n)]=0;~~~\sum_{k<n}\xi_k{\rm tr }[c(w)\widehat{c}(\widetilde{e}_k)\widehat{c}(\widetilde{e}_n)c(\xi')]=0.
\end{align} Then, we have
\begin{align}\label{8888}
&{\rm tr} [\pi^+_{\xi_n}\sigma_{0}(-2c(w)\nabla^{\bigwedge^*T^*M} _v\widetilde{D}^{-1})\times
\frac{i}{2(1+\xi_n^2)^2}
  h'(0)\sum_{k<n}\xi_kc(\widetilde{e}_k)c(\widetilde{e}_n)](x_0)\nonumber\\
&=\frac{i}{2(\xi_n-i)^4(\xi_n+i)^2}h'(0)\sum_{j=1}^{n-1}\sum_{k<n}v_j\xi_j\xi_kw_k{\rm tr}[\texttt{id}]-\frac{i}{2(\xi_n-i)^4(\xi_n+i)^2}h'(0)\sum_{j=1}^{n-1}\sum_{k<n}\xi_j\xi_kv_nw_n{\rm tr}[\texttt{id}],
\end{align}
\begin{align}\label{hhh888}
{\rm tr} [\pi^+_{\xi_n}\sigma_{0}(-2c(w)\nabla^{\bigwedge^*T^*M} v\widetilde{D}^{-1})\times
\frac{i}{2(1+\xi_n^2)^2}
 h'(0)\sum_{k<n}\xi_k\widehat{c}(\widetilde{e}_k)\widehat{c}(\widetilde{e}_n)](x_0)=0,
\end{align}
and 
\begin{align}\label{8kk8}
&{\rm tr} [\pi^+_{\xi_n}\sigma_{0}(-2c(w)\nabla^{\bigwedge^*T^*M} v\widetilde{D}^{-1})\times
-\frac{5i\xi_n^3+9i\xi_n}{2(1+\xi_n^2)^3}h'(0)](x_0)\nonumber\\
&=-\frac{5\xi_n^3+9\xi_n}{2(\xi_n-i)^5(\xi_n+i)^3}h'(0)\sum_{j,k=1}^{n-1}v_jw_k\xi_j\xi_k{\rm tr}[\texttt{id}]+\frac{5\xi_n^3+9\xi_n}{2(\xi_n-i)^5(\xi_n+i)^3}h'(0)v_nw_n{\rm tr}[\texttt{id}].
\end{align}
Next, we perform the corresponding integral calculation on the above results. Therefore
\begin{align}\label{39}
\Phi_4^b&=i\int_{|\xi'|=1}\int^{+\infty}_{-\infty}{\rm tr} [\partial_{\xi_n}\pi^+_{\xi_n}\sigma_{0}(-2c(w)\nabla^{\bigwedge^*T^*M} _v\widetilde{D}^{-1})\times
\sigma_{-3}(\widetilde{D}^{-2})](x_0)d\xi_n\sigma(\xi')dx'\nonumber\\
&=\frac{2\pi i}{3!}\left[\frac{1}{2(\xi_n+i)^2}\right]^{(3)}\bigg|_{\xi_n=i}\times -\frac{4\pi}{3}\sum_{j=1}^{n-1}h'(0)\times 16dx'+\frac{2\pi i}{3!}\left[\frac{1}{ 2(\xi_n+i)^2}\right]^{(3)}\bigg|_{\xi_n=i}\times v_nw_nh'(0)\times 16\Omega_3dx'\nonumber\\
&\frac{2\pi i}{4!}\left[\frac{5i\xi_n^3+9i\xi_n}{2(\xi_n+i)^3}\right]^{(4)}\bigg|_{\xi_n=i}\times  -\frac{4\pi}{3}\sum_{j=1}^{n-1}h'(0)\times 16dx'-\frac{2\pi i}{4!}\left[\frac{5i\xi_n^3+9i\xi_n}{2(\xi_n+i)^3}\right]^{(4)}\bigg|_{\xi_n=i}\times -v_nw_n\times 16\Omega_3dx'\nonumber\\
 &=\bigg(\frac{38}{3}g(v^T,w^T)-\frac{98}{3}v_nw_n\bigg)h'(0)\pi^2dx'.
\end{align}
\noindent {\bf  case c)}~$r=-1,~\ell=-2,~k=j=|\alpha|=0$.\\
By (\ref{a3}), we get
\begin{align}\label{61}
\Phi^b_5=-i\int_{|\xi'|=1}\int^{+\infty}_{-\infty}{\rm tr} [\pi^+_{\xi_n}\sigma_{-1}(-2c(w)\nabla^{\bigwedge^*T^*M} _v\widetilde{D}^{-1}
)\times
\partial_{\xi_n}\sigma_{-2}(\widetilde{D}^{-2})](x_0)d\xi_n\sigma(\xi')dx'.
\end{align}
By Lemma \ref{lemb1}, we have
\begin{equation}\label{62}
\partial_{\xi_n}\sigma_{-2}(\widetilde{D}^{-2})(x_0)|_{|\xi'|=1}=\frac{-2\xi_n}{(1+\xi_n^2)^2}.
\end{equation}
By Lemma \ref{lemb2e}, we have
\begin{align}\label{63}
\sigma_{-1}(-2c(w)\nabla^{\bigwedge^*T^*M} _v\widetilde{D}^{-1})(x_0):=A_1(x_0)+A_2(x_0)+A_3(x_0),
\end{align}
where
\begin{align}
A_1(x_0)&=-2c(w)\sqrt{-1}\sum_{j=1}^nv_j\xi_j\bigg[\frac{c(\xi)\sigma_{0}(\widetilde{D})c(\xi)}{|\xi|^4}+\frac{c(\xi)}{|\xi|^6}\sum_jc(dx_j)
\Big(\partial_{x_j}(c(\xi))|\xi|^2-c(\xi)\partial_{x_j}(|\xi|^2)\Big)\bigg];\nonumber\\
A_2(x_0)&=-2c(w)A(v)\frac{\sqrt{-1}c(\xi)}{|\xi|^2};\nonumber\\
A_3(x_0)&=-2c(w)v_n\bigg(\frac{\sqrt{-1}\partial_{x_n}c(\xi')}{|\xi|^2}-\frac{\sqrt{-1}c(\xi)|\xi'|^2h'(0)}{|\xi|^4}\bigg)(x_0),
\end{align}
where
\begin{align}\label{4ss4}
\sigma_{0}(\widetilde{D})(x_0)&=\frac{1}{4}\sum_{s,t,i}\omega_{s,t}(\widetilde{e}_i)
(x_{0})c(\widetilde{e}_i)\widehat{c}(\widetilde{e}_s)\widehat{c}(\widetilde{e}_t)
-\frac{1}{4}\sum_{s,t,i}\omega_{s,t}(\widetilde{e}_i)
(x_{0})c(\widetilde{e}_i)c(\widetilde{e}_s)c(\widetilde{e}_t).
\end{align}
We denote
\begin{align}\label{45}
Q_0^{1}(x_0)&=\frac{1}{4}\sum_{s,t,i}\omega_{s,t}(\widetilde{e}_i)
(x_{0})c(\widetilde{e}_i)\widehat{c}(\widetilde{e}_s)\widehat{c}(\widetilde{e}_t);\nonumber\\
Q_0^{2}(x_0)&=-\frac{1}{4}\sum_{s,t,i}\omega_{s,t}(\widetilde{e}_i)
(x_{0})c(\widetilde{e}_i)c(\widetilde{e}_s)c(\widetilde{e}_t).
\end{align}
Firstly, the following results are obtained by further calculation of $A_1(x_0$)
\begin{align}
A_1(x_0)&=\sum_{j=1}^{n-1}(-2c(w)v_j\sqrt{-1}\xi_j)\bigg(\frac{3\xi_n^4+4\xi_n^2-7}{4(1+\xi_n^2)^3}h'(0)c(dx_n)+\frac{3\xi_n^2+7\xi_n}{2(1+\xi_n^2)^3}h'(0)c(\xi')\nonumber\\
&+\frac{1}{(1+\xi_n^2)^2}c(\xi')c(dx_n)\partial_{x_n}(c(\xi'))-\frac{\xi_n}{(1+\xi_n^2)^2}\bigg)-2c(w)v_n\sqrt{-1}\bigg(\frac{3\xi_n^5+4\xi_n^3-7\xi_n}{4(1+\xi_n^2)^3}h'(0)c(dx_n)\nonumber\\
&+\frac{3\xi_n^3+7\xi_n^2}{2(1+\xi_n^2)^3}h'(0)c(\xi')+\frac{\xi_n}{(1+\xi_n^2)^2}c(\xi')c(dx_n)\partial_{x_n}(c(\xi'))-\frac{\xi_n^2}{(1+\xi_n^2)^2}\bigg).
\end{align}
If we omit some items that have no contribution for computing $\Phi^b_5$, by the Cauchy integral formula, we obtain
\begin{align}\label{64}
\pi^+_{\xi_n}A_1(x_0)&=\frac{i}{(\xi_n-i)^3}\sum_{j=1}^{n-1}v_j\xi_jh'(0)c(w)c(\xi')-\frac{i}{2(\xi_n-i)^3}v_nh'(0)c(w)c(dx_n)\nonumber\\
&-\frac{1}{2(\xi_n-i)^2}v_n\xi_j\sum_{j=1}^{n-1}v_j\xi_jc(w)\partial_{x_n}(c(\xi'))+\frac{1}{2(\xi_n-i)^3}v_nh'(0)c(w)c(\xi')c(dx_n)\partial_{x_n}(c(\xi')).
\end{align}
Since
\begin{align}\label{4bbb8}
c(dx_n)Q_0^{1}(x_0)
&=-\frac{1}{4}h'(0)\sum^{n-1}_{i=1}c(\widetilde{e}_i)
\widehat{c}(\widetilde{e}_i)c(\widetilde{e}_n)\widehat{c}(\widetilde{e}_n);\nonumber\\
 Q_0^{2}&=c_0c(dx_n)=-\frac{3}{4}h'(0)c(dx_n).
\end{align}
By the relation of the Clifford action, we have the following identities:
\begin{align}\label{ggg32}
&\sum_{j=1}^{n-1}\xi_jv_j{\rm tr }[c(w)\partial_{x_n}(c(\xi'))]=-\frac{h'(0)}{2}\sum_{j,k=1}^{n-1}\xi_j\xi_kv_jw_k{\rm tr}[\texttt{id}];~~{\rm tr}[c(w)c(\xi')c(dx_n)\partial_{x_n}(c(\xi'))]=\frac{h'(0)}{2}w_n{\rm tr}[\texttt{id}];\nonumber\\
&\sum_{j=1}^{n-1}v_j\xi_j{\rm tr}[c(w)c(\xi')Q_0^{1}c(dx_n)](x_0)|_{|\xi'|=1}=0;~~~{\rm tr}[c(w)c(\xi')Q_0^{1}c(\xi')](x_0)|_{|\xi'|=1}=0;\nonumber\\
&{\rm tr}[c(w)c(dx_n)Q_0^{1}c(dx_n)](x_0)|_{|\xi'|=1}=0.
\end{align}
Then
 \begin{align}\label{65}
&-i\int_{|\xi'|=1}\int^{+\infty}_{-\infty}{\rm tr} [\pi^+_{\xi_n}A_1(x_0)\times
\partial_{\xi_n}\sigma_{-2}(\widetilde{D}^{-2})](x_0)d\xi_n\sigma(\xi')dx'\nonumber\\
&=\bigg(-\frac{16}{3}g(v^T,w^T)+2v_nw_n\bigg)h'(0)\pi^2dx'.
\end{align}
Secondly, for $A_2(x_0)$, further calculation leads to new results
\begin{align}\label{66}
\pi^+_{\xi_n} A_2(x_0)&=\pi^+_{\xi_n}\Big(-2c(w)A(v)\frac{\sqrt{-1}c(\xi)}{|\xi|^2}\Big)\nonumber\\
&=-\frac{1}{\xi_n-i}c(w)A(v)c(\xi')-\frac{i}{\xi_n-i}c(w)A(v)c(dx_n).
\end{align}
Next
\begin{align}
&{\rm tr} [\pi^+_{\xi_n}A_2(x_0)\times
\partial_{\xi_n}\sigma_{-2}(\widetilde{D}^{-2})](x_0)\nonumber\\
&=\frac{2\xi_n}{(\xi_n-i)(1+\xi_n^2)^2}{\rm tr}[c(w)A(v)c(\xi')]+\frac{2i\xi_n}{(\xi_n-i)(1+\xi_n^2)^2}{\rm tr}[c(w)A(v)c(dx_n)].
\end{align}
When $i<n,~\int_{|\xi'|=1}\xi_{i_{1}}\xi_{i_{2}}\cdots\xi_{i_{2d+1}}\sigma(\xi')=0$, ${\rm tr }[c(w)A(v)c(\xi')]$ has no contribution for computing $\Phi^b_5$.
Then we only need to compute  ${\rm tr}[c(w)A(v)c(dx_n)]$, by the relation of the Clifford action, we have the following identities
\begin{align*}
&{\rm tr}[c(w)A(v)c(dx_n)]\nonumber\\
&=\frac{1}{4}\sum_{i,j=1}^n{\rm tr}[c(w)<\nabla^L_v\widetilde{e}_i,\widetilde{e}_j>c(\widetilde{e}_i)c(\widetilde{e}_j)c(dx_n)]-\frac{1}{4}\sum_{i,j=1}^n{\rm tr}[c(w)<\nabla^L_v\widetilde{e}_i,\widetilde{e}_j>\widehat{c}(\widetilde{e}_i)\widehat{c}(\widetilde{e}_j)c(dx_n)]\nonumber\\
&=\frac{1}{4}\sum_{i,j=1}^n{\rm tr}[<\nabla^L_v\widetilde{e}_i,\widetilde{e}_j>c(\widetilde{e}_i)c(\widetilde{e}_j)c(dx_n)c(w)]-\frac{1}{4}\sum_{i,j=1}^n{\rm tr}[<\nabla^L_v\widetilde{e}_i,\widetilde{e}_j>\widehat{c}(\widetilde{e}_i)\widehat{c}(\widetilde{e}_j)c(dx_n)c(w)]\nonumber\\
&=\frac{1}{4}\sum_{1\leq i,j\leq n-1}{\rm tr}[<\nabla^L_v\widetilde{e}_i,\widetilde{e}_j>c(\widetilde{e}_i)c(\widetilde{e}_j)c(dx_n)c(w)]-\frac{1}{4}\sum_{1\leq i,j\leq n-1}{\rm tr}[<\nabla^L_v\widetilde{e}_i,\widetilde{e}_j>\widehat{c}(\widetilde{e}_i)\widehat{c}(\widetilde{e}_j)c(dx_n)c(w)]\nonumber\\
&+\frac{1}{4}\sum_{j\leq n-1}{\rm tr}[<\nabla^L_v\frac{\partial}{\partial{x_n}},\widetilde{e}_j>c(dx_n)c(\widetilde{e}_j)c(dx_n)c(w)]-\frac{1}{4}\sum_{j\leq n-1}{\rm tr}[<\nabla^L_v\frac{\partial}{\partial{x_n}},\widetilde{e}_j>\widehat{c}(dx_n)\widehat{c}(\widetilde{e}_j)c(dx_n)c(w)]\nonumber\\
&+\frac{1}{4}\sum_{i\leq n-1}{\rm tr}[<\nabla^L_v\widetilde{e}_i,\frac{\partial}{\partial{x_n}}>c(\widetilde{e}_i)c(dx_n)c(dx_n)c(w)]-\frac{1}{4}\sum_{i\leq n-1}{\rm tr}[<\nabla^L_v\widetilde{e}_i,\frac{\partial}{\partial{x_n}}>\widehat{c}(\widetilde{e}_i)\widehat{c}(dx_n)c(dx_n)c(w)]\nonumber\\
&+\frac{1}{4}{\rm tr}[<\nabla^L_v\frac{\partial}{\partial{x_n}},\frac{\partial}{\partial{x_n}}>c(dx_n)c(dx_n)c(dx_n)c(w)]-\frac{1}{4}{\rm tr}[<\nabla^L_v\frac{\partial}{\partial{x_n}},\frac{\partial}{\partial{x_n}}>\widehat{c}(dx_n)\widehat{c}(dx_n)c(dx_n)c(w)].
\end{align*}
By $<\nabla^L_v\widetilde{e}_i,\widetilde{e}_j>+<\widetilde{e}_i,\nabla^L_v\widetilde{e}_j>=v<\widetilde{e}_i,\widetilde{e}_j>$, we have\\
(1)when $i=j,$ 
\begin{align*}
<\nabla^L_v\widetilde{e}_i,\widetilde{e}_j>=0;
\end{align*}
(2)when $i\neq j\leq n-1,$ 
\begin{align*}
&\sum_{i\neq j\leq n-1}{\rm tr}[<\nabla^L_v\widetilde{e}_i,\widetilde{e}_j>c(\widetilde{e}_i)c(\widetilde{e}_j)c(dx_n)c(w)]\nonumber\\
&=\sum_{i\neq j\leq n-1}\sum_{l=1}^{n-1}{\rm tr}[<\nabla^L_v\widetilde{e}_i,\widetilde{e}_j>c(\widetilde{e}_i)c(\widetilde{e}_j)c(dx_n)w_lc(\widetilde{e}_l)]+\sum_{i\neq j\leq n-1}{\rm tr}[<\nabla^L_v\widetilde{e}_i,\widetilde{e}_j>c(\widetilde{e}_i)c(\widetilde{e}_j)c(dx_n)w_nc(\widetilde{e}_n)]\nonumber\\
&=0;
\end{align*}
\begin{align*}
&-\sum_{i\neq j\leq n-1}{\rm tr}[<\nabla^L_v\widetilde{e}_i,\widetilde{e}_j>\widehat{c}(\widetilde{e}_i)\widehat{c}(\widetilde{e}_j)c(dx_n)c(w)]\nonumber\\
&=-\sum_{i\neq j\leq n-1}\sum_{l=1}^{n-1}{\rm tr}[<\nabla^L_v\widetilde{e}_i,\widetilde{e}_j>\widehat{c}(\widetilde{e}_i)\widehat{c}(\widetilde{e}_j)c(dx_n)w_lc(\widetilde{e}_l)]-\sum_{i\neq j\leq n-1}{\rm tr}[<\nabla^L_v\widetilde{e}_i,\widetilde{e}_j>\widehat{c}(\widetilde{e}_i)\widehat{c}(\widetilde{e}_j)c(dx_n)w_nc(\widetilde{e}_n)]\nonumber\\
&=0;
\end{align*}
(3)when $=n, j\leq n-1,$
\begin{align*}
\sum_{j\leq n-1}{\rm tr}[<\nabla^L_v\frac{\partial}{\partial{x_n}},\widetilde{e}_j>c(\widetilde{e}_j)c(dx_n)c(dx_n)c(w)]=-\sum_{j\leq n-1}<\nabla^L_v\frac{\partial}{\partial{x_n}},\widetilde{e}_j>w_j{\rm tr}[\texttt{id}];
\end{align*}
\begin{align*}
-\sum_{j\leq n-1}{\rm tr}[<\nabla^L_v\frac{\partial}{\partial{x_n}},\widetilde{e}_j>\widehat{c}(\widetilde{e}_j)\widehat{c}(dx_n)c(dx_n)c(w)]=0.
\end{align*}
Then
\begin{align}\label{1a}
{\rm tr}[c(w)A(v)c(dx_n)]&=-\frac{1}{2}\sum_{j\leq n-1}<\nabla^L_v\frac{\partial}{\partial{x_n}},\widetilde{e}_j>w_j{\rm tr}[\texttt{id}]=-\frac{1}{2}<\nabla^L_v\frac{\partial}{\partial{x_n}},w^T>{\rm tr}[\texttt{id}].
\end{align}
Moreover
\begin{align}
&-i\int_{|\xi'|=1}\int^{+\infty}_{-\infty}{\rm tr} [\pi^+_{\xi_n}A_2(x_0)\times
\partial_{\xi_n}\sigma_{-2}(\widetilde{D}^{-2})](x_0)d\xi_n\sigma(\xi')dx'\nonumber\\
&=-i\int_{\Gamma^+}\frac{i\xi_n}{(\xi_n-i)^3(\xi_n+i)^2}d\xi_n<\nabla^L_v\frac{\partial}{\partial{x_n}},w^T>{\rm tr}[\texttt{id}]dx'\nonumber\\
&=-i\times\frac{2\pi i}{2!}\left[\frac{i\xi_n}{(\xi_n+i)^2}\right]^{(2)}\bigg|_{\xi_n=i}<\nabla^L_v\frac{\partial}{\partial{x_n}},w^T>\times16\times4\pi dx'\nonumber\\
&=8<\nabla^L_v\frac{\partial}{\partial{x_n}},w^T>\pi^2dx'.
\end{align}
Thirdly, for $A_3(x_0)$, we get
\begin{align}\label{66pp}
\pi^+_{\xi_n} A_3(x_0)&=\pi^+_{\xi_n}\bigg[-2c(w)v_n\bigg(\frac{\sqrt{-1}\partial_{x_n}c(\xi')}{|\xi|^2}-\frac{\sqrt{-1}c(\xi)|\xi'|^2h'(0)}{|\xi|^4}\bigg)\bigg]\nonumber\\
&=-\frac{1}{\xi_n-i}v_nc(w)\partial_{x_n}[c(\xi')]+\frac{1}{2(\xi_n-i)}h'(0)v_nc(w)c(dx_n).
\end{align}
By (\ref{66pp}) and ${\rm tr }[c(w)c(\xi')]$ has no contribution for computing $\Phi^b_5$, we have
\begin{align}
&-i\int_{|\xi'|=1}\int^{+\infty}_{-\infty}{\rm tr} [\pi^+_{\xi_n}A_3(x_0)\times
\partial_{\xi_n}\sigma_{-2}(\widetilde{D}^{-2})](x_0)d\xi_n\sigma(\xi')dx'\nonumber\\
&=-i\int_{|\xi'|=1}\int^{+\infty}_{-\infty} \frac{\xi_n}{(\xi_n-i)^4(\xi_n+i)^2}h'(0)v_nw_n{\rm tr}[\texttt{id}]d\xi_n\sigma(\xi')dx'\nonumber\\
&=-i\int_{\Gamma^+}\frac{\xi_n}{(\xi_n-i)^4(\xi_n+i)^2}h'(0)v_nw_n{\rm tr}[\texttt{id}]d\xi_n\Omega_3dx'\nonumber\\
&=-i\frac{2\pi i}{3!}\left[\frac{\xi_n}{(\xi_n+i)^2}\right]^{(3)}\bigg|_{\xi_n=i}h'(0)v_nw_n\times16\times4\pi dx'\nonumber\\
&=8h'(0)v_nw_n\pi^2dx'.
\end{align}
Therefore
\begin{align}\label{6666}
\Phi^b_5&=-i\int_{|\xi'|=1}\int^{+\infty}_{-\infty}{\rm tr} [\pi^+_{\xi_n}(A_1+A_2+A_3)(x_0)\times
\partial_{\xi_n}\sigma_{-2}(\widetilde{D}^{-2})](x_0)d\xi_n\sigma(\xi')dx'\nonumber\\
&=\bigg(-\frac{16}{3}h'(0)g(v^T,w^T)+8<\nabla^L_v\frac{\partial}{\partial{x_n}},w^T>+10h'(0)v_nw_n\bigg)\pi^2dx'.
\end{align}
Because $\Phi^b$ is the sum of the cases (a), (b) and (c). Finally, we get
\begin{align}\label{1795}
\Phi^b&=\bigg(-\frac{8}{3}\partial_{x_n}g(v^T,w^T)+8\partial_{x_n}(v_nw_n)-\frac{44}{3}h'(0)v_nw_n+\frac{18}{3}h'(0)g(v^T,w^T)+8<\nabla^L_v\frac{\partial}{\partial{x_n}},w^T> \bigg)\pi^2dx'.
\end{align}
By (3.78) in \cite{WJ2}, we know that when $n=4,$ the following identity holds
\begin{align}\label{ppop}
K(x_0)=\sum_{ij}K_{ij}(x_0)g_{\partial M}^{ij}(x_0)=\sum_{i=1}^3K_{ii}(x_0)=-\frac{3}{2}h'(0),
\end{align}
where
$K_{ij}$ is the second fundamental form, or extrinsic
curvature.

Substituting (\ref{ppop}) into (\ref{1795}), we have
\begin{align}\label{17195}
\Phi^b&=\bigg(-\frac{8}{3}\partial_{x_n}g(v^T,w^T)+8\partial_{x_n}(v_nw_n)+\frac{88}{9}Kv_nw_n-\frac{36}{9}Kg(v^T,w^T)+8<\nabla^L_v\frac{\partial}{\partial{x_n}},w^T>\bigg)\pi^2dx'.
\end{align}
Combine the results of boundary $\Phi^a$ and boundary $\Phi^b$, we obtain following theorem
\begin{thm}\label{thmb1}
Let M be a $4$-dimensional compact oriented Riemannian manifold with boundary $\partial M$ and the metric
$g^M$ be defined in Section \ref{section:2}, then we get the following equality:
\begin{align}
\label{b2773}
&\widetilde{{\rm Wres}}[\pi^+(c(w)(\widetilde{D}c(v)+c(v)\widetilde{D})\widetilde{D}^{-1})\circ\pi^+(\widetilde{D}^{-2})]\nonumber\\
&=\frac{64\pi^2}{3}\int_M[Ric(v,w)-\frac{1}{2}s(g)g(v,w)]{Vol_g}+\int_{\partial M}\bigg\{-8\bigg(\sum_{j=1}^ng(e_j,\nabla^L_{e_j}v)g(w,\frac{\partial}{\partial{x_n}})-g(w,\nabla^L_{\frac{\partial}{\partial{x_n}}}v)\nonumber\\
&+g(\nabla^L_wv,\frac{\partial}{\partial{x_n}})\bigg)-\frac{8}{3}\partial_{x_n}g(v^T,w^T)+8\partial_{x_n}(v_nw_n)+\frac{88}{9}Kv_nw_n-\frac{36}{9}Kg(v^T,w^T)\nonumber\\
&+8<\nabla^L_v\frac{\partial}{\partial{x_n}},w^T>\bigg\}\pi^2dx'.
\end{align}
\end{thm}
\section{ The type-II operator and Dabrowski-Sitarz-Zalecki type theorems for 4-dimensional manifolds with boundary}
\label{section:5}
In this section, we give the the Einstein functional about the type-II operator for $4$ dimensional manifold with boundary, and prove the Dabrowski-Sitarz-Zalecki type theorem about the type-II operator.
 
By Propsition 4.1 in \cite{WJ2}, we get the following propsition
\begin{prop}\label{1propb1}For the type-II operator, the Einstein functional for $4$ dimensional spin manifolds with boundary defined by
\begin{align}\label{1b1}
&\widetilde{{\rm Wres}}[\pi^+(c(w)(\widetilde{D}c(v)+c(v)\widetilde{D})\widetilde{D}^{-2})\circ\pi^+(\widetilde{D}^{-1})]\nonumber\\
&=\widetilde{{\rm Wres}}[\pi^+(\sum_{j=1}^nc(w)c(\widetilde{e}_j)c(\nabla^{T^*M}_{\widetilde{e}_j}v)\widetilde{D}^{-2})\circ\pi^+(\widetilde{D}^{-1})]+\widetilde{{\rm Wres}}[\pi^+(-2c(w)\nabla^{\bigwedge^*T^*M} _v\widetilde{D}^{-2})\circ\pi^+(\widetilde{D}^{-1})].
\end{align}
\end{prop}
\subsection{$\widetilde{{\rm Wres}}[\pi^+(\sum_{j=1}^nc(w)c(\widetilde{e}_j)c(\nabla^{T^*M}_{\widetilde{e}_j}v)\widetilde{D}^{-2})\circ\pi^+(\widetilde{D}^{-1})]$}
Combining with the generating Proposition \ref{propb1}, this yields
\begin{align}\label{1b2}
\widetilde{{\rm Wres}}[\pi^+(\sum_{j=1}^nc(w)c(\widetilde{e}_j)c(\nabla^{T^*M}_{\widetilde{e}_j}v)\widetilde{D}^{-2})\circ\pi^+(\widetilde{D}^{-1})]={\rm Wres}[\sum_{j=1}^nc(w)c(\widetilde{e}_j)c(\nabla^{T^*M}_{\widetilde{e}_j}v)\widetilde{D}^{-3}]+\int_{\partial M}\Psi^a,\nonumber\\
\end{align}

Now, we can compute $\Psi^a$. The sum is taken over $
r+l-k-j-|\alpha|=-3,~~r=-2,~~l=-1,$ then we have the only one case:
\begin{align}\label{1b4}
\Psi^a=-i\int_{|\xi'|=1}\int^{+\infty}_{-\infty}
 {\rm tr}[\pi^+_{\xi_n}\sigma_{-2}(\sum_{j=1}^nc(w)c(\widetilde{e}_j)c(\nabla^{T^*M}_{\widetilde{e}_j}v)\widetilde{D}^{-2})\times\partial_{\xi_n}\sigma_{-1}(\widetilde{D}^{-1})](x_0)d\xi_n\sigma(\xi')dx'.
\end{align}
By Lemma \ref{lemb1}, we get
\begin{align}
\pi^+_{\xi_n}\sigma_{-2}(\sum_{j=1}^nc(w)c(\widetilde{e}_j)c(\nabla^{T^*M}_{\widetilde{e}_j}v)\widetilde{D}^{-2})=\frac{-i}{2(\xi_n-i)}\sum_{j=1}^nc(w)c(\widetilde{e}_j)c(\nabla^{T^*M}_{\widetilde{e}_j}v),\nonumber
\end{align}
and by further derivative, we get
\begin{align}
\partial_{\xi_n}\sigma_{-1}(\widetilde{D}^{-1})=\frac{-2i\xi_n}{(1+\xi_n^2)^2}c(\xi')+\frac{i(1-\xi_n^2)}{(1+\xi_n^2)^2}c(dx_n).\nonumber
\end{align}
It follows that 
\begin{align}
&{\rm tr}[\pi^+_{\xi_n}\sigma_{-2}(\sum_{j=1}^nc(w)c(\widetilde{e}_j)c(\nabla^{T^*M}_{\widetilde{e}_j}v)\widetilde{D}^{-2})\times\partial_{\xi_n}\sigma_{-1}(\widetilde{D}^{-1})](x_0)\nonumber\\
&=\frac{-\xi_n}{(\xi_n-i)^2(1+\xi_n^2)^2}{\rm tr}\sum_{j=1}^n[c(w)c(\widetilde{e}_j)c(\nabla^{T^*M}_{\widetilde{e}_j}v)c(\xi')](x_0)\nonumber\\
&+\frac{1-\xi_n^2}{2(\xi_n-i)(1+\xi_n^2)^2}{\rm tr}\sum_{j=1}^n[c(w)c(\widetilde{e}_j)c(\nabla^{T^*M}_{\widetilde{e}_j}v)c(dx_n)](x_0)\nonumber.
\end{align}
By (\ref{123}) and ${\rm tr}[c(w)c(e_j)c(\nabla^{T*M}_{\widetilde{e}_j}v)c(\xi')](x_0)$ has no contribution for computing $\Psi^a$. Then, we have
\begin{align}\label{1iiii}
\Psi^a&=-i\int_{|\xi'|=1}\int^{+\infty}_{-\infty}
 {\rm tr}[\pi^+_{\xi_n}\sigma_{-2}(\sum_{j=1}^nc(w)c(\widetilde{e}_j)c(\nabla^{T^*M}_{\widetilde{e}_j}v)\widetilde{D}^{-2})\times\partial_{\xi_n}\sigma_{-1}(\widetilde{D}^{-1})](x_0)d\xi_n\sigma(\xi')dx'\nonumber\\
&=\int_{\Gamma^+}\frac{i\xi_n^2-i}{2(\xi_n-i)^3(\xi_n+i)^2}d\xi_n\bigg(\sum_{j=1}^ng(e_j,\nabla^L_{e_j}v)g(w,\frac{\partial}{\partial{x_n}})-g(w,\nabla^L_{\frac{\partial}{\partial{x_n}}}v)+g(\nabla^L_wv,\frac{\partial}{\partial{x_n}})\bigg){\rm {\rm tr}[\texttt{id}]}\Omega_3dx'\nonumber\\
&=\frac{2\pi i}{2!}\left[\frac{i\xi_n^2-i}{(\xi_n+i)^2}\right]^{(2)}\bigg(\sum_{j=1}^ng(e_j,\nabla^L_{e_j}v)g(w,\frac{\partial}{\partial{x_n}})-g(w,\nabla^L_{\frac{\partial}{\partial{x_n}}}v)+g(\nabla^L_wv,\frac{\partial}{\partial{x_n}})\bigg)\times16\times4\pi dx'\nonumber\\
&=8\bigg(\sum_{j=1}^ng(e_j,\nabla^L_{e_j}v)g(w,\frac{\partial}{\partial{x_n}})-g(w,\nabla^L_{\frac{\partial}{\partial{x_n}}}v)+g(\nabla^L_wv,\frac{\partial}{\partial{x_n}})\bigg)\pi^2dx'.
\end{align}

\subsection{$\widetilde{{\rm Wres}}[\pi^+(-2c(w)\nabla^{\bigwedge^*T^*M}
_v\widetilde{D}^{-2})\circ\pi^+(\widetilde{D}^{-1})]$}
Combining with the generating Proposition \ref{propb1}, this yields
\begin{align}
\widetilde{{\rm Wres}}[\pi^+(-2c(w)\nabla^{\bigwedge^*T^*M}
_v\widetilde{D}^{-2})\circ\pi^+(\widetilde{D}^{-1})]={\rm Wres}[-2c(w)\nabla^{\bigwedge^*T^*M}
_v\widetilde{D}^{-3}]+\int_{\partial M}\Psi^b.\nonumber
\end{align}

By Lemma \ref{lemb1}, Lemma \ref{lem553} and the composition formula of pseudodifferential operators, we have the following lemma
\begin{lem}\label{1lemb2e}The following identities hold:
\begin{align}
\sigma_{-1}(\nabla_v^{\bigwedge^*T^*M}
\widetilde{D}^{-2})&=\sqrt{-1}\sum_{j=1}^nv_j\xi_j|\xi|^{-2};\nonumber\\
\sigma_{-2}(\nabla_v^{\bigwedge^*T^*M}
\widetilde{D}^{-2})&=\sigma_0(\nabla_v^{\bigwedge^*T^*M}
)\sigma_{-2}(\widetilde{D}^{-2})+\sigma_{1}(\nabla_v^{\bigwedge^*T^*M}
)\sigma_{-3}(\widetilde{D}^{-2})+\sum_{j=1}^n\xi_j\sigma_{1}(\nabla_v^{S(TM)})D_{x_j}\sigma_{-2}(\widetilde{D}^{-2})\nonumber\\
&=A(v)|\xi|^{-2}-\sqrt{-1}\sum_{j=1}^nv_j\xi_j\bigg[\sqrt{-1}|\xi|^4\xi_k(\Gamma^k-2\sigma^k+2a^k)+2\sqrt{-1}\xi|^{-6}\xi^j\xi_\alpha\xi_\beta\partial_jg^{\alpha\beta}\bigg]\nonumber\\
&-v_n\frac{h'(0)|\xi'|^2}{|\xi|^4}.\nonumber
\end{align}
\end{lem}

 The same to (\ref{llll}), we have
\begin{align}\label{llrrll}
{\rm Wres}[c(w)(\widetilde{D}c(v)+c(v)\widetilde{D})\widetilde{D}^{-3}]=\frac{64\pi^2}{3}\int_M[Ric(v,w)-\frac{1}{2}s(g)g(v,w)]{Vol_g}.
\end{align}
 
Therefore, we only need to compute $\Psi^b$. The sum is taken over $
r+l-k-j-|\alpha|=-3,~~r\leq -1,~~l\leq -1,$ then we have the following five cases:
~\\
\noindent  {\bf case a)~I)}~$r=-1,~l=-1,~k=j=0,~|\alpha|=1$.\\
\noindent By (\ref{a5}), we get
\begin{align}
\Psi^b_1=-\int_{|\xi'|=1}\int^{+\infty}_{-\infty}\sum_{|\alpha|=1}
 {\rm tr}[\partial^\alpha_{\xi'}\pi^+_{\xi_n}\sigma_{-1}(-2c(w)\nabla^{\bigwedge^*T^*M}
_v\widetilde{D}^{-2})\times
 \partial^\alpha_{x'}\partial_{\xi_n}\sigma_{-1}(\widetilde{D}^{-1})](x_0)d\xi_n\sigma(\xi')dx'.\nonumber
\end{align}
By Lemma 2.2 in \cite{Wa3}, for $i<n$, then
\begin{align}
\partial_{x_i}\sigma_{-1}(\widetilde{D}^{-1})(x_0)=\partial_{x_i}\frac{ic(\xi')}{|\xi|^2}(x_0)=0,\nonumber
\end{align}
so $\Psi^b_1=0$.\\
 \noindent  {\bf case a)~II)}~$r=-1,~l=-1,~k=|\alpha|=0,~j=1$.\\
\noindent By (\ref{a5}), we get
\begin{align}
\Psi^b_2=-\frac{1}{2}\int_{|\xi'|=1}\int^{+\infty}_{-\infty} {\rm
tr} [\partial_{x_n}\pi^+_{\xi_n}\sigma_{-1}(-2c(w)\nabla^{\bigwedge^*T^*M}
_v\widetilde{D}^{-2})\times
\partial_{\xi_n}^2\sigma_{-1}(\widetilde{D}^{-1})](x_0)d\xi_n\sigma(\xi')dx'.
\end{align}
\noindent By Lemma \ref{lemb1}, we have
\begin{align}
\partial_{\xi_n}^2\sigma_{-1}(\widetilde{D}^{-1})(x_0)=\partial_{\xi_n}^2\left(\frac{ic(\xi')}{|\xi|^{2}}\right)(x_0)=\frac{i(6\xi_n^2-2)}{(1+\xi_n^2)^3}c(\xi')+\frac{i(2\xi_n^3-6\xi_n)}{(1+\xi_n^2)^3}c(dx_n).\nonumber
\end{align}
By Lemma \ref{1lemb2e}, we have\\
\begin{align}\label{1b27}
&\partial_{x_n}\sigma_{-1}(-2c(w)\nabla_v^{S(TM)}\widetilde{D}^{-2})\nonumber\\
&=-2i\sum_{j=1}^{n-1}\xi_j\bigg[\frac{\partial_{x_n}(v_jc(w))}{|\xi|^2}-\frac{h'(0)|\xi'|^2}{|\xi|^4}v_jc(w)\bigg]-2i\xi_n\bigg[\frac{\partial_{x_n}(v_nc(w))}{|\xi|^2}-\frac{h'(0)|\xi'|^2}{|\xi|^4}v_nc(w)\bigg].
\end{align}
Moreover
\begin{align}\label{1b28}
&\partial_{x_n}\pi^+_{\xi_n}\sigma_{-1}(-2c(w)\nabla_v^{\bigwedge^*T^*M}
\widetilde{D}^{-2})\nonumber\\
&=\pi^+_{\xi_n}\partial_{x_n}\sigma_{-1}(-2c(w)\nabla_v^{\bigwedge^*T^*M}
\widetilde{D}^{-2})\nonumber\\
&=\sum_{j=1}^{n-1}\xi_j\bigg[\frac{-1}{(\xi_n-i)}\partial_{x_n}(v_jc(w))-\frac{2+i\xi_n}{2(\xi_n-i)^2}h'(0)v_jc(w)\bigg]-\frac{i}{\xi_n-i}\partial_{x_n}(v_nc(w))+\frac{1}{2(\xi_n-i)^2}h'(0)v_nc(w).
\end{align}
By (\ref{eee}) and we omit some items that have no contribution for computing $\Psi^b_2$. Then, we have
\begin{align}\label{1b78}
& {\rm
tr} [\partial_{x_n}\pi^+_{\xi_n}\sigma_{-1}(-2c(w)\nabla^{\bigwedge^*T^*M}_v\widetilde{D}^{-2})\times
\partial_{\xi_n}^2\sigma_{-1}(\widetilde{D}^{-1})](x_0)\nonumber\\
&=-\frac{i(2-6\xi_n^2)}{(\xi_n-i)(1+\xi_n^2)^3}\sum_{j,k=1}^{n-1}\xi_j\xi_k\partial_{x_n}(v_jw_k){\rm tr}[\texttt{id}]-\frac{(2+i\xi_n)(3\xi_n^2-1)}{(\xi_n-i)^2(1+\xi_n^2)^3}h'(0)|\xi'|^2\sum_{j,k=1}^{n-1}\xi_j\xi_kv_jw_k{\rm tr}[\texttt{id}]\nonumber\\
&-\frac{2\xi_n^3-6\xi_n}{(\xi_n-i)(1+\xi_n^2)^3}\partial_{x_n}(v_nw_n){\rm tr}[\texttt{id}]-\frac{i\xi_n^3-3i\xi_n}{(\xi_n-i)^2(1+\xi_n^2)^3}h'(0)|\xi'|^2v_nw_n{\rm tr}[\texttt{id}].
\end{align}
Therefore, we have
\begin{align}\label{135}
\Psi^b_2&=-\frac{1}{2}\int_{|\xi'|}\int_{-\infty}^{\infty}{\rm
tr} [\partial_{x_n}\pi^+_{\xi_n}\sigma_{-1}(-2c(w)\nabla^{S(TM)}_v\widetilde{D}^{-2})\times
\partial_{\xi_n}^2\sigma_{-1}(\widetilde{D}^{-1})](x_0)d\xi_n\sigma(\xi')dx'\nonumber\\
&=-\frac{1}{2}\times\frac{2\pi i}{3!}\left[\frac{i(6\xi_n-2)}{(\xi_n+i)^3}\right]^{(3)}\bigg|_{\xi_n=i}\times \frac{4\pi}{3}\sum_{j=1}^{n-1}\partial_{x_n}(v_jw_j)h'(0)\times 16dx'\nonumber\\
&+\frac{1}{2}\times\frac{2\pi i}{4!}\left[\frac{3i\xi_n^3+6\xi_n^2-i\xi_n-2}{(\xi_n+i)^3}\right]^{(4)}\bigg|_{\xi_n=i}\times \frac{4\pi}{3}\sum_{j=1}^{n-1}v_jw_j\times 16dx'\nonumber\\
&-\frac{1}{2}\times\frac{2\pi i}{3!}\left[\frac{6\xi_n-2\xi_n^3}{(\xi_n+i)^3}\right]^{(3)}\bigg|_{\xi_n=i}\times \partial_{x_n}(v_nw_n)\times 16\Omega_3dx'\nonumber\\
&-\frac{1}{2}\times\frac{2\pi i}{4!}\left[\frac{3i\xi_n-\xi_n^3}{(\xi_n+i)^3}\right]^{(4)}\bigg|_{\xi_n=i}\times v_nw_n\times 16\Omega_3dx'\nonumber\\
 &=\bigg(-\frac{8}{3}\partial_{x_n}g(v^T,w^T)+\frac{10}{3}h'(0)g(v^T,w^T)-8\partial_{x_n}(v_nw_n)+6h'(0)v_nw_n\bigg)\pi^2dx'.
 \end{align}
\noindent  {\bf case a)~III)}~$r=-1,~l=-1,~j=|\alpha|=0,~k=1$.\\
\noindent By (\ref{a5}), we get
\begin{align}\label{136}
\Psi^b_3&=-\frac{1}{2}\int_{|\xi'|=1}\int^{+\infty}_{-\infty}
{\rm tr} [\partial_{\xi_n}\pi^+_{\xi_n}\sigma_{-1}(-2c(w)\nabla^{\bigwedge^*T^*M}
_v\widetilde{D}^{-2})\times
\partial_{\xi_n}\partial_{x_n}\sigma_{-1}(\widetilde{D}^{-1})](x_0)d\xi_n\sigma(\xi')dx'.
\end{align}
By Lemma \ref{lemb1} and by further calculation, we have
\begin{align}\label{138887}
\partial_{\xi_n}\partial_{x_n} \sigma_{-1}(D^{-1})(x_0)|_{|\xi'|=1}=\frac{4i\xi_nh'(0)}{(1+\xi_n^2)^3}c(\xi')+\frac{3i\xi_n^2-i}{(1+\xi_n^2)^3}c(dx_n)-\frac{2i\xi_n}{(1+\xi_n^2)^3}\partial_{x_n}(c(\xi')).
\end{align}
By (\ref{64})
Then, we get
\begin{align}\label{1mmmmm}
\partial_{\xi_n}\pi^+_{\xi_n}\sigma_{-1}(-2c(w)\nabla^{\bigwedge^*T^*M}
_v\widetilde{D}^{-2})
&=\sum_{j=1}^{n-1}v_j\xi_j\frac{1}{(\xi_n-i)^2}c(w)-v_n\frac{i}{(\xi_n-i)^2}c(w).
\end{align}
By (\ref{eee}), we have
\begin{align}\label{1kkk}
&{\rm tr} [\partial_{\xi_n}\pi^+_{\xi_n}\sigma_{-1}(-2c(w)\nabla^{\bigwedge^*T^*M}
_v\widetilde{D}^{-2})\times
\partial_{\xi_n}\partial_{x_n}\sigma_{-1}(\widetilde{D}^{-1})](x_0)\nonumber\\
&=-\frac{4i\xi_n}{(\xi_n-i)^2(1+\xi_n^2)^3}h'(0)\sum_{j,k=1}^{n-1}v_jw_k\xi_j\xi_k{\rm tr}[\texttt{id}]+\frac{i\xi_n}{(\xi_n-i)^2(1+\xi_n^2)^2}h'(0)\sum_{j,k=1}^{n-1}v_jw_k\xi_j\xi_k{\rm tr}[\texttt{id}]\nonumber\\
&-\frac{(1-3\xi_n^2)}{(\xi_n-i)^2(1+\xi_n^2)^3}h'(0)v_nw_n{\rm tr}[\texttt{id}].
\end{align}
Next, we perform the corresponding integral calculation on the above results. Therefore
\begin{align}\label{13ll5}
\Psi^b_3&=-\frac{1}{2}\int_{|\xi'|=1}\int^{+\infty}_{-\infty}
{\rm tr} [\partial_{\xi_n}\pi^+_{\xi_n}\sigma_{0}(-2c(w)\nabla^{\bigwedge^*T^*M}
_v\widetilde{D}^{-1})\times
\partial_{\xi_n}\partial_{x_n}\sigma_{-2}(\widetilde{D}^{-2})](x_0)d\xi_n\sigma(\xi')dx'\nonumber\\
&=-\frac{1}{2}\times\frac{2\pi i}{4!}\left[\frac{4i\xi_n}{(\xi_n+i)^3}\right]^{(4)}\bigg|_{\xi_n=i}\times-\frac{4\pi}{3}\sum_{j=1}^{n-1}v_jw_j\times 16h'(0)dx'\nonumber\\
&-\frac{1}{2}\times\frac{2\pi i}{4!}\left[\frac{3\xi_n^2-1}{(\xi_n+i)^3}\right]^{(4)}\bigg|_{\xi_n=i}\times v_nw_n\times 16h'(0)\Omega_3dx'\nonumber\\
 &=\bigg(-2g(v^T,w^T)-6v_nw_n\bigg)h'(0)\pi^2dx'.
 \end{align}
\noindent  {\bf case b)}~$r=-2,~l=-1,~k=j=|\alpha|=0$.\\
\noindent By (\ref{a5}), we get
\begin{align}\label{142}
\Psi^b_4&=-i\int_{|\xi'|=1}\int^{+\infty}_{-\infty}{\rm tr} [\pi^+_{\xi_n}\sigma_{-2}(-2c(w)\nabla^{\bigwedge^*T^*M}
_v\widetilde{D}^{-2})\times
\partial_{\xi_n}\sigma_{-1}(\widetilde{D}^{-1})](x_0)d\xi_n\sigma(\xi')dx'.
\end{align}
By Lemma \ref{lemb1}, we have
\begin{align}\label{143}
\partial_{\xi_n}\sigma_{-1}(\widetilde{D}^{-1})](x_0)|_{|\xi'|=1}=\sqrt{-1}\bigg[\frac{c(dx_n)}{1+\xi_n^2}-\frac{2\xi_nc(\xi')+2\xi_n^2c(dx_n)}{(1+\xi_n^2)^2}\bigg].
\end{align}
By Lemma \ref{1lemb2e}, we have
\begin{align}\label{673}
\sigma_{-2}(-2c(w)\nabla^{S(TM)}_v\widetilde{D}^{-2})(x_0):=B_1(x_0)+B_2(x_0)+B_3(x_0),
\end{align}
where
\begin{align}
B_1(x_0)&=-2c(w)A(v)|\xi|^{-2};\nonumber\\
B_2(x_0)&=-2\sum_{j=1}^n\sqrt{-1}v_j\xi_j\bigg[\frac{i}{2(1+\xi_n^2)^2}
   h'(0)\sum_{k<n}\xi_kc(\widetilde{e}_k)c(\widetilde{e}_n)\nonumber\\
   &+\frac{i}{2(1+\xi_n^2)^2}h'(0)\sum_{k<n}\xi_k\widehat{c}(\widetilde{e}_k)\widehat{c}(\widetilde{e}_n)-\frac{5i\xi_n^3+9i\xi_n}{2(1+\xi_n^2)^3}h'(0)\bigg];\nonumber\\
B_3(x_0)&=2c(w)v_n\frac{h'(0)|\xi'|^2}{|\xi|^4}.
\end{align}
Firstly, the following results are obtained by further calculation of $B_1(x_0)$
\begin{align}
\pi^+_{\xi_n} B_1(x_0)=\frac{i}{\xi_n-i}c(w)A(v).
\end{align}
Then
\begin{align}
{\rm tr} [\pi^+_{\xi_n} B_1\times
\partial_{\xi_n}\sigma_{-1}(\widetilde{D}^{-1})](x_0)=-\frac{1-\xi_n^2}{(\xi_n-i)(1+\xi_n^2)^2}{\rm tr}[c(w)A(v)c(dx_n)]+\frac{\xi_n}{(\xi_n-i)(1+\xi_n^2)^2}{\rm tr}[c(w)A(v)c(\xi')].
\end{align}
We note that ${\rm tr}[c(w)A(v)c(\xi')]$ no contribution for computing $\Psi^b_4$, then
 \begin{align}\label{65}
&-i\int_{|\xi'|=1}\int^{+\infty}_{-\infty}{\rm tr} [\pi^+_{\xi_n}B_1(x_0)\times
\partial_{\xi_n}\sigma_{-2}(\widetilde{D}^{-2})](x_0)d\xi_n\sigma(\xi')dx'\nonumber\\
&=i\int_{\Gamma^+}\frac{1-\xi_n^2}{(\xi_n-i)^3(\xi_n+i)^2}d\xi_n\times-\frac{1}{2}<\nabla^L_v\frac{\partial}{\partial{x_n}},w^T> {\rm tr}[\texttt{id}]\Omega_3dx'\nonumber\\
&=i\times\frac{2\pi i}{2!}\left[\frac{1-\xi_n^2}{(\xi_n+i)^2}\right]^{(2)}\bigg|_{\xi_n=i}\times-\frac{1}{2}<\nabla^L_v\frac{\partial}{\partial{x_n}},w^T> \times16\times4\pi dx'\nonumber\\
&=-8<\nabla^L_v\frac{\partial}{\partial{x_n}},w^T> \pi^2dx'.
\end{align}
Secondly, for $B_2(x_0)$, further calculation leads to new results
\begin{align}\label{66}
\pi^+_{\xi_n} B_2(x_0)&=-\frac{2+i\xi_n}{4(\xi_n-i)^2}h'(0)\sum_{j=1}^{n-1}v_j\xi_j\sum_{k<n}\xi_kc(w)c(\widetilde{e}_k)c(\widetilde{e}_n)-\frac{i}{4(\xi_n-i)^2}h'(0)v_n\sum_{k<n}\xi_kc(w)c(\widetilde{e}_k)c(\widetilde{e}_n)\nonumber\\
&-\frac{2+i\xi_n}{4(\xi_n-i)^2}h'(0)\sum_{j=1}^{n-1}v_j\xi_j\sum_{k<n}\xi_kc(w)\widehat{c}(\widetilde{e}_k)\widehat{c}(\widetilde{e}_n)-\frac{i}{4(\xi_n-i)^2}h'(0)v_n\sum_{k<n}\xi_kc(w)\widehat{c}(\widetilde{e}_k)\widehat{c}(\widetilde{e}_n)\nonumber\\
&+\frac{1}{2(\xi_n-i)^3}h'(0)\sum_{j=1}^{n-1}v_j\xi_jc(w)+\frac{i}{2(\xi_n-i)^3}h'(0)v_nc(w).
\end{align}
Moreover
\begin{align}
&{\rm tr} [\pi^+_{\xi_n}B_2(x_0)\times
\partial_{\xi_n}\sigma_{-2}(\widetilde{D}^{-2})](x_0)\nonumber\\
&=\frac{-(2+i\xi_n)(i-i\xi_n^2)}{4(\xi_n-i)^2(1+\xi_n^2)^2}h'(0)\sum_{j=1,k<n}^{n-1}v_j\xi_j\xi_k{\rm tr} [c(w)c(\widetilde{e}_k)c(\widetilde{e}_n)c(dx_n)]+\frac{(\xi_n^2-1)h'(0)}{2(\xi_n-i)^3(1+\xi_n^2)^2}v_n{\rm tr} [c(w)c(dx_n)]\nonumber\\
&-\frac{\xi_n}{2(\xi_n-i)^2(1+\xi_n^2)^2}h'(0)v_n\sum_{k<n}\xi_k{\rm tr} [c(w)c(\widetilde{e}_k)c(\widetilde{e}_n)c(\xi')]-\frac{i\xi_n}{(\xi_n-i)^3(1+\xi_n^2)^2}h'(0)\sum_{j=1}^{n-1}v_j\xi_j{\rm tr} [c(w)c(\xi')]\nonumber\\
&\frac{-(2+i\xi_n)(i-i\xi_n^2)}{4(\xi_n-i)^2(1+\xi_n^2)^2}h'(0)\sum_{j=1,k<n}^{n-1}v_j\xi_j\xi_k{\rm tr} [c(w)\widehat{c}(\widetilde{e}_k)\widehat{c}(\widetilde{e}_n)c(dx_n)]\nonumber\\
&-\frac{\xi_n}{2(\xi_n-i)^2(1+\xi_n^2)^2}h'(0)v_n\sum_{k<n}\xi_k{\rm tr} [c(w)\widehat{c}(\widetilde{e}_k)\widehat{c}(\widetilde{e}_n)c(\xi')].
\end{align}
Then by (\ref{lll32}), we have
\begin{align}
&-i\int_{|\xi'|=1}\int^{+\infty}_{-\infty}{\rm tr} [\pi^+_{\xi_n}B_2(x_0)\times
\partial_{\xi_n}\sigma_{-2}(\widetilde{D}^{-2})](x_0)d\xi_n\sigma(\xi')dx'\nonumber\\
&=-i\int_{\Gamma^+}\frac{-(2+i\xi_n)(i-i\xi_n^2)}{4(\xi_n-i)^2(1+\xi_n^2)^2}h'(0)d\xi_n\int_{|\xi'|=1}\sum_{j,k=1}^{n-1}v_jw_k\xi_j\xi_k\sigma(\xi'){\rm tr}[\texttt{id}]dx'\nonumber\\
&-i\int_{\Gamma^+}\frac{(\xi_n^2-1)h'(0)}{2(\xi_n-i)^3(1+\xi_n^2)^2}h'(0)d\xi_n\int_{|\xi'|=1}\times-v_nw_n{\rm tr}[\texttt{id}]dx'\nonumber\\
&-i\int_{\Gamma^+}\frac{\xi_n}{2(\xi_n-i)^2(1+\xi_n^2)^2}h'(0)d\xi_n\int_{|\xi'|=1}\sum_{j,k=1}^{n-1}v_nw_n\xi_j\xi_k\sigma(\xi'){\rm tr}[\texttt{id}]dx'\nonumber\\
&-i\int_{\Gamma^+}\frac{i\xi_n}{(\xi_n-i)^3(1+\xi_n^2)^2}h'(0)d\xi_n\int_{|\xi'|=1}\sum_{j,k=1}^{n-1}v_jw_k\xi_j\xi_k\sigma(\xi'){\rm tr}[\texttt{id}]dx'\nonumber\\
&=\frac{2}{3}h'(0)\sum_{j=1}^{n-1}v_jw_j\pi^2dx'-\frac{22}{3}h'(0)v_nw_n\pi^2dx'.
\end{align}
Thirdly, for $B_3(x_0)$, we get
\begin{align}\label{661pp}
\pi^+_{\xi_n} B_3(x_0)
&=-\frac{2+i\xi_n}{2(\xi_n-i)^2}h'(0)v_nc(w).
\end{align}
Then
\begin{align}
&{\rm tr} [\pi^+_{\xi_n}B_3(x_0)\times
\partial_{\xi_n}\sigma_{-2}(\widetilde{D}^{-2})](x_0)\nonumber\\
&=\frac{2i\xi_n-\xi_n^2}{2(\xi_n-i)^2(1+\xi_n^2)^2}h'(0)v_n{\rm tr}[c(w)c(\xi')]-\frac{\xi_n^3-2i\xi_n^2-\xi_n+2i}{2(\xi_n-i)^2(1+\xi_n^2)^2}v_n{\rm tr}[c(w)c(dx_n)].
\end{align}
By (\ref{66pp}) and ${\rm tr }[c(w)c(\xi')]$ has no contribution for computing $\Psi^b_4$, we have
\begin{align}
&-i\int_{|\xi'|=1}\int^{+\infty}_{-\infty}{\rm tr} [\pi^+_{\xi_n}B_3(x_0)\times
\partial_{\xi_n}\sigma_{-2}(\widetilde{D}^{-2})](x_0)d\xi_n\sigma(\xi')dx'\nonumber\\
&=-i\int_{|\xi'|=1}\int^{+\infty}_{-\infty} \frac{\xi_n^3-2i\xi_n^2-\xi_n+2i}{2(\xi_n-i)^2(1+\xi_n^2)^2}h'(0)v_nw_n{\rm tr}[\texttt{id}]d\xi_n\sigma(\xi')dx'\nonumber\\
&=-i\int_{\Gamma^+} \frac{\xi_n^3-2i\xi_n^2-\xi_n+2i}{2(\xi_n-i)^2(1+\xi_n^2)^2}h'(0)v_nw_n{\rm tr}[\texttt{id}]d\xi_n\Omega_3dx'\nonumber\\
&=-i\frac{2\pi i}{3!}\left[\frac{\xi_n^3-2i\xi_n^2-\xi_n+2i}{2(\xi_n+i)^2}\right]^{(3)}\bigg|_{\xi_n=i}h'(0)v_nw_n\times16\times4\pi dx'\nonumber\\
&=-16h'(0)v_nw_n\pi^2dx'.
\end{align}
Therefore
\begin{align}\label{61666}
\Psi^b_4&=--i\int_{|\xi'|=1}\int^{+\infty}_{-\infty}{\rm tr} [\pi^+_{\xi_n}(B_1+B_2+B_3)(x_0)\times
\partial_{\xi_n}\sigma_{-2}(\widetilde{D}^{-2})](x_0)d\xi_n\sigma(\xi')dx'\nonumber\\
&=\bigg(\frac{2}{3}h'(0)g(v^T,w^T)-\frac{70}{3}h'(0)v_nw_n-8<\nabla^L_v\frac{\partial}{\partial{x_n}},w^T> \bigg)\pi^2dx'.
\end{align}
\noindent {\bf  case c)}~$r=-1,~\ell=-2,~k=j=|\alpha|=0$.\\
By (\ref{a5}), we get
\begin{align}\label{161}
\Psi^b_5=-i\int_{|\xi'|=1}\int^{+\infty}_{-\infty}{\rm tr} [\pi^+_{\xi_n}\sigma_{-1}(-2c(w)\nabla^{\bigwedge^*T^*M}
_v\widetilde{D}^{-2}
)\times
\partial_{\xi_n}\sigma_{-2}(\widetilde{D}^{-1})](x_0)d\xi_n\sigma(\xi')dx'.
\end{align}
 By Lemma \ref{1lemb2e}, we have
 \begin{align}
 \sigma_{-1}(-2c(w)\nabla^{\bigwedge^*T^*M}
_v\widetilde{D}^{-2})=-2\sqrt{-1}c(w)\left(\sum_{j=1}^{n-1}v_j\xi_j|\xi|^{-2}+v_n\xi_n|\xi|^{-2}\right)
 \end{align}
 By the Cauchy integral formula, we obtain
  \begin{align}
\pi^+_{\xi_n}\sigma_{-1}(-2c(w)\nabla^{\bigwedge^*T^*M}
_v\widetilde{D}^{-2})=-\frac{1}{\xi_n-i}\sum_{j=1}^{n-1}v_j\xi_jc(w)-\frac{i}{\xi_n-i}v_nc(w).
 \end{align}
 By Lemma \ref{lemb1}, we get
 \begin{align}\label{4ss3}
\sigma_{-2}(\widetilde{D}^{-1})(x_0)=\frac{c(\xi)\sigma_{0}(\widetilde{D})(x_0)c(\xi)}{|\xi|^4}+\frac{c(\xi)}{|\xi|^6}c(dx_n)
[\partial_{x_n}(c(\xi'))(x_0)|\xi|^2-c(\xi)h'(0)|\xi|^2_{\partial
M}]
\end{align}
Then by (\ref{4ss4}), we have
\begin{align}\label{65}
&\partial_{\xi_n}\sigma_{-2}(\widetilde{D}^{-1})(x_0)|_{|\xi'|=1}\nonumber\\
&=
\partial_{\xi_n}\bigg\{\frac{c(\xi)[Q_0^{2}(x_0)+Q_0^{2}(x_0)]c(\xi)}{|\xi|^4}+\frac{c(\xi)}{|\xi|^6}c(dx_n)[\partial_{x_n}[c(\xi')](x_0)|\xi|^2-c(\xi)h'(0)]\bigg\}\nonumber\\
&=\partial_{\xi_n}\bigg\{\frac{[c(\xi)Q_0^{1}(x_0)]c(\xi)}{|\xi|^4}+\frac{c(\xi)}{|\xi|^6}c(dx_n)[\partial_{x_n}[c(\xi')](x_0)|\xi|^2-c(\xi)h'(0)]\bigg\}\nonumber\\
&+\partial_{\xi_n}\frac{c(\xi)Q_0^{2}(x_0)c(\xi)}{|\xi|^4}\nonumber\\
&:=C_1+C_2,
\end{align}
where 
\begin{align}
C_1&=\partial_{\xi_n}\bigg\{\frac{[c(\xi)Q_0^{1}(x_0)]c(\xi)}{|\xi|^4}+\frac{c(\xi)}{|\xi|^6}c(dx_n)[\partial_{x_n}[c(\xi')](x_0)|\xi|^2-c(\xi)h'(0)]\bigg\};\nonumber\\
C_2&=\partial_{\xi_n}\frac{c(\xi)Q_0^{2}(x_0)c(\xi)}{|\xi|^4}.
\end{align}
Firstly, for $C_1$, further calculation leads to new results
\begin{align}
C_1&=\frac{1}{(1+\xi_n^2)^3}\bigg[(2\xi_n-2\xi_n^3)c(dx_n)Q_0^{2}c(dx_n)
+(1-3\xi_n^2)c(dx_n)Q_0^{2}c(\xi')+(1-3\xi_n^2)c(\xi')Q_0^{2}c(dx_n)\nonumber\\
&-4\xi_nc(\xi')Q_0^{2}c(\xi')
+(3\xi_n^2-1){\partial}_{x_n}c(\xi')-4\xi_nc(\xi')c(dx_n){\partial}_{x_n}c(\xi')+2h'(0)c(\xi')+2h'(0)\xi_nc(dx_n)\bigg]\nonumber\\
&+6\xi_nh'(0)\frac{c(\xi)c(dx_n)c(\xi)}{(1+\xi^2_n)^4}.
\end{align}
Moreover
\begin{align}
&{\rm tr} [\pi^+_{\xi_n}\sigma_{-1}(-2c(w)\nabla^{\bigwedge^*T^*M}
_v\widetilde{D}^{-2}
)\times
C_1](x_0)\nonumber\\
&=\frac{3(1-3\xi_n^2)}{2(\xi_n-i)(1+\xi_n^2)^3}h'(0)\sum_{j,k=1}^{n-1}v_jw_k\xi_j\xi_k{\rm tr}[\texttt{id}]+\frac{1-3\xi_n^2}{2(\xi_n-i)(1+\xi_n^2)^3}h'(0)\sum_{j,k=1}^{n-1}v_jw_k\xi_j\xi_k{\rm tr}[\texttt{id}]\nonumber\\
&+\frac{2}{(\xi_n-i)(1+\xi_n^2)^3}h'(0)\sum_{j,k=1}^{n-1}v_jw_k\xi_j\xi_k{\rm tr}[\texttt{id}]-\frac{12\xi_n^2}{(\xi_n-i)(1+\xi_n^2)^4}h'(0)\sum_{j,k=1}^{n-1}v_jw_k\xi_j\xi_k{\rm tr}[\texttt{id}]\nonumber\\
&+\frac{3i(\xi_n-\xi_n^3)}{2(\xi_n-i)(1+\xi_n^2)^3}h'(0)v_nw_n{\rm tr}[\texttt{id}]+\frac{3i\xi_n}{(\xi_n-i)(1+\xi_n^2)^3}h'(0)v_nw_n{\rm tr}[\texttt{id}]\nonumber\\
&-\frac{2i\xi_n}{(\xi_n-i)(1+\xi_n^2)^3}h'(0)v_nw_n{\rm tr}[\texttt{id}]+\frac{2i\xi_n}{(\xi_n-i)(1+\xi_n^2)^3}h'(0)v_nw_n{\rm tr}[\texttt{id}]\nonumber\\
&+\frac{6i\xi_n}{(\xi_n-i)(1+\xi_n^2)^4}h'(0)v_nw_n{\rm tr}[\texttt{id}]-\frac{6i\xi_n^3}{(\xi_n-i)(1+\xi_n^2)^4}h'(0)v_nw_n{\rm tr}[\texttt{id}].
\end{align}
Then
\begin{align}
&-i\int_{|\xi'|=1}\int^{+\infty}_{-\infty}{\rm tr} [\pi^+_{\xi_n}\sigma_{-1}(-2c(w)\nabla^{\bigwedge^*T^*M}
_v\widetilde{D}^{-2}
)\times
C_1](x_0)d\xi_n\sigma(\xi')dx'\nonumber\\
&=\bigg(6\sum_{j=1}^{n-1}v_jw_j+18v_nw_n\bigg)h'(0)\pi^2dx'.
\end{align}
Secondly, by further calculation of $C_2$, the following results are obtained
\begin{align}
C_2&=\frac{c(dx_n)Q_0^{1}(x_0)c(\xi)}{|\xi|^4}
+\frac{c(\xi)Q_0^{1}(x_0)c(dx_n)}{|\xi|^4}
-\frac{4\xi_n c(\xi)Q_0^{1}(x_0)c(\xi)}{|\xi|^6}.
\end{align}
By (\ref{ggg32}), and we omit some items that have no contribution for computing $\Psi^b_5$, we get
\begin{align}
&{\rm tr} [\pi^+_{\xi_n}\sigma_{-1}(-2c(w)\nabla^{\bigwedge^*T^*M}
_v\widetilde{D}^{-2}
)\times
C_2](x_0)\nonumber\\
&=-\frac{1-3\xi_n^2}{(\xi_n-i)(1+\xi_n^2)^3}\sum_{j=1}^{n-1}v_j\xi_j{\rm tr} [c(w)c(dx_n)Q_0^{1}(x_0)c(\xi')]-\frac{i(1-3\xi_n^2)}{(\xi_n-i)(1+\xi_n^2)^3}v_n{\rm tr} [c(w)c(dx_n)Q_0^{1}(x_0)c(\xi')]\nonumber\\
&-\frac{2(\xi_n-\xi_n^3)}{(\xi_n-i)(1+\xi_n^2)^3}\sum_{j=1}^{n-1}v_j\xi_j{\rm tr} [c(w)c(dx_n)Q_0^{1}(x_0)c(dx_n)]-\frac{2i(\xi_n-\xi_n^3)}{(\xi_n-i)(1+\xi_n^2)^3}v_n{\rm tr} [c(w)c(dx_n)Q_0^{1}(x_0)c(dx_n)]\nonumber\\
&-\frac{1-3\xi_n^2}{(\xi_n-i)(1+\xi_n^2)^3}\sum_{j=1}^{n-1}v_j\xi_j{\rm tr} [c(w)c(\xi')Q_0^{1}(x_0)c(dx_n)]-\frac{i(1-3\xi_n^2)}{(\xi_n-i)(1+\xi_n^2)^3}v_n{\rm tr} [c(w)c(\xi')Q_0^{1}(x_0)c(dx_n)]\nonumber\\
&+\frac{4\xi_n}{(\xi_n-i)(1+\xi_n^2)^3}\sum_{j=1}^{n-1}v_j\xi_j{\rm tr} [c(w)c(\xi')Q_0^{1}(x_0)c(\xi')]+\frac{4i\xi_n}{(\xi_n-i)(1+\xi_n^2)^3}v_n{\rm tr} [c(w)c(\xi')Q_0^{1}(x_0)c(\xi')]\nonumber\\
&=0.
\end{align}
Then
\begin{align}
-i\int_{|\xi'|=1}\int^{+\infty}_{-\infty}{\rm tr} [\pi^+_{\xi_n}\sigma_{-1}(-2c(w)\nabla^{\bigwedge^*T^*M}
_v\widetilde{D}^{-2}
)\times
C_2](x_0)d\xi_n\sigma(\xi')dx'=0.
\end{align}
Therefore
\begin{align}\label{61666}
\Psi^b_5&=-i\int_{|\xi'|=1}\int^{+\infty}_{-\infty}{\rm tr} [\pi^+_{\xi_n}\sigma_{-1}(-2c(w)\nabla^{\bigwedge^*T^*M}
_v\widetilde{D}^{-2}
)\times
(C_1+C_2)](x_0)d\xi_n\sigma(\xi')dx'\nonumber\\
&=\bigg(6g(v^T,w^T)+18v_nw_n\bigg)h'(0)\pi^2dx'.
\end{align}
Now $\Psi^b$ is the sum of the cases (a), (b) and (c). Therefore, we get
\begin{align}\label{795}
\Psi^b&=\bigg(-\frac{8}{3}\partial_{x_n}g(v^T,w^T)-8\partial_{x_n}(v_nw_n)-\frac{16}{3}h'(0)v_nw_n+8h'(0)g(v^T,w^T)-8<\nabla^L_v\frac{\partial}{\partial{x_n}},w^T> \bigg)\pi^2dx'.
\end{align}
By $K(x_0)=-\frac{3}{2}h'(0)$, we have
\begin{align}\label{7935}
\Psi^b&=\bigg(-\frac{8}{3}\partial_{x_n}g(v^T,w^T)-8\partial_{x_n}(v_nw_n)+\frac{32}{9}Kv_nw_n-\frac{16}{9}Kg(v^T,w^T)-8<\nabla^L_v\frac{\partial}{\partial{x_n}},w^T>\bigg)\pi^2dx'.
\end{align}
Combine the results of boundary $\Psi^a$ and boundary $\Psi^b$, we obtain following theorem
\begin{thm}\label{1thmb1}
Let M be a $4$-dimensional compact oriented Riemannian manifold with boundary $\partial M$ and the metric
$g^M$ be defined as above, then we get the following equality:
\begin{align}
\label{b27173}
&\widetilde{{\rm Wres}}[\pi^+c(w)(\widetilde{D}c(v)+c(v)\widetilde{D})\widetilde{D}^{-2}\circ\pi^+(\widetilde{D}^{-1})]\nonumber\\
&=\frac{64\pi^2}{3}\int_M[Ric(v,w)-\frac{1}{2}s(g)g(v,w)]{Vol_g}+\int_{\partial M}\bigg\{\bigg(8\sum_{j=1}^ng(e_j,\nabla^L_{e_j}v)g(w,\frac{\partial}{\partial{x_n}})-g(w,\nabla^L_{\frac{\partial}{\partial{x_n}}}v)\nonumber\\
&+g(\nabla^L_wv,\frac{\partial}{\partial{x_n}})\bigg)-\frac{8}{3}\partial_{x_n}g(v^T,w^T)-8\partial_{x_n}(v_nw_n)+\frac{32}{9}Kv_nw_n-\frac{16}{9}Kg(v^T,w^T)\nonumber\\
&-8<\nabla^L_v\frac{\partial}{\partial{x_n}},w^T>\bigg\}\pi^2dx'.
\end{align}
\end{thm}
\section*{Acknowledgements}
This work was supported by NSFC. 11771070. The authors thank the referee for his (or her) careful reading and helpful comments.

\section*{}

\end{document}